\documentclass{birkjour}

\usepackage[utf8]{inputenc}
\usepackage[T1]{fontenc}
\usepackage[english]{babel}
\usepackage{mathtools,amsmath,amssymb,amsfonts,amsthm}
\usepackage{bbm}
\mathtoolsset{showonlyrefs}
\usepackage{enumitem}

\usepackage{geometry}
\usepackage{graphicx,color}
\usepackage{tikz,pgfplots,pgf}
\usepackage{epsfig}
\usepackage{subcaption}
\usepackage[font=small]{caption}

\definecolor{c1}{RGB}{89,126,183}
\definecolor{c2}{RGB}{235,155,0}
\definecolor{c3}{RGB}{127,184,45}
\definecolor{c4}{RGB}{255,81,10}
\definecolor{c5}{RGB}{142,112,180}

\newtheorem{theorem}{Theorem}[section]

\newtheorem{proposition}{Proposition}[section]
\newtheorem{corollary}{Corollary}[section]

\newtheorem{conjecture}{Conjecture}[section]
\theoremstyle{definition}
\newtheorem{remark}{Remark}[section]

\numberwithin{equation}{section}

\begin{document}

\title[Uniqueness of positive solutions for indefinite boundary value problems]{Uniqueness of positive solutions \\for boundary value problems associated with \\indefinite $\phi$-Laplacian type equations}

\author[A.~Boscaggin]{Alberto Boscaggin}

\address{
Department of Mathematics ``Giuseppe Peano'', University of Torino\\
Via Carlo Alberto 10, 10123 Torino, Italy}

\email{alberto.boscaggin@unito.it}

\author[G.~Feltrin]{Guglielmo Feltrin}

\address{
Department of Mathematics, Computer Science and Physics, University of Udine\\
Via delle Scienze 206, 33100 Udine, Italy}

\email{guglielmo.feltrin@uniud.it}

\author[F.~Zanolin]{Fabio Zanolin}

\address{
Department of Mathematics, Computer Science and Physics, University of Udine\\
Via delle Scienze 206, 33100 Udine, Italy}

\email{fabio.zanolin@uniud.it}

\thanks{Work written under the auspices of the Grup\-po Na\-zio\-na\-le per l'Anali\-si Ma\-te\-ma\-ti\-ca, la Pro\-ba\-bi\-li\-t\`{a} e le lo\-ro Appli\-ca\-zio\-ni (GNAMPA) of the Isti\-tu\-to Na\-zio\-na\-le di Al\-ta Ma\-te\-ma\-ti\-ca (INdAM). The first two authors are supported by INdAM--GNAMPA project ``Problemi ai limiti per l'equazione della curvatura media prescritta''.
\\
\textbf{Preprint -- September 2020}}

\subjclass{34B15, 34B16, 34B18, 34C25.}

\keywords{Uniqueness, indefinite weight, positive solutions, $p$-Laplacian, boundary value problems, superlinear functions, singular equations.}

\date{}

\dedicatory{}

\begin{abstract}
The paper provides a uniqueness result for positive solutions of the Neumann and periodic boundary value problems associated with the $\phi$-Laplacian equation
\begin{equation*}
\bigl{(} \phi(u') \bigr{)}' + a(t) g(u) = 0,
\end{equation*}
where $\phi$ is a homeomorphism with $\phi(0)=0$, $a(t)$ is a stepwise indefinite weight and $g(u)$ is a continuous function. When dealing with the $p$-Laplacian differential operator $\phi(s)=|s|^{p-2}s$ with $p>1$, and the nonlinear term $g(u)=u^{\gamma}$ with $\gamma\in\mathbb{R}$, we prove the existence of a unique positive solution when $\gamma\in\mathopen{]}-\infty,(1-2p)/(p-1)\mathclose{]} \cup \mathopen{]}p-1,+\infty\mathclose{[}$.
\end{abstract}

\maketitle

\section{Introduction}\label{section-1}

This paper deals with the $\phi$-Laplacian differential equation
\begin{equation}\label{eq-phi}
\bigl{(}\phi(u') \bigr{)}' + a(t) g(u) = 0,
\end{equation}
where $\phi$ is an increasing homeomorphism defined in an open interval including the origin, with $\phi(0)=0$, $a(t)$ is a sign-changing $L^{1}$-weight function and $g(u)$ is a continuous function with $g(u)>0$ for $u>0$.

It is worth noticing that the $\phi$-Laplacian operator appearing in equation~\eqref{eq-phi} includes several classical differential operators, such as the linear operator $\phi(s)=s$, or the $p$-Laplacian operator $\phi(s)=\varphi_{p}(s)=|s|^{p-2}s$ with $p>1$, or even the one-dimensional mean-curvature operator in Lorentz--Minkowski spaces $\phi(s)=s/\sqrt{1-|s|^{2}}$. Such differential operators have been widely investigated in the literature for their relevance in many mechanical and physical models (cf.~\cite{KrRaVa-10}).

Starting from \cite{HeKa-80}, it is common to refer to \eqref{eq-phi} as a nonlinear indefinite equation, due to the presence of a sign-changing weight function. The study of indefinite problems, both in the ODE and in the PDE setting, has shown an exceptional interest, from the pioneering works \cite{AlTa-96,AmLG-98,BeCDNi-95} till to the recent developments dealing with positive solutions of boundary value problems associated with \eqref{eq-phi} (we refer to \cite{Fe-18book} for a quite comprehensive list of references).

In this paper, we focus our attention on the Neumann and the periodic boundary value problems associated with \eqref{eq-phi} and we analyse both the power-type nonlinearity, that is
\begin{equation}\label{g-super}
g(u)=u^{q}, \quad q>0,
\end{equation}
and the singular nonlinearity
\begin{equation}\label{g-sing}
g(u)=\dfrac{1}{u^{\kappa}}, \quad \kappa>0.
\end{equation}
In this framework, lot of work has been done concerning existence and multiplicity of positive solutions, see, for instance, \cite{BaPoTe-88,BeCDNi-95,Bo-11,BoFe-20jde,FeZa-15ade,FeZa-17jde,LGOm-20,LGOmRi-17jde} for the power-type nonlinearity \eqref{g-super}, and \cite{BeZa-18,BoZa-15,BrTo-10,GoZa-19jdde,GoZa-19,HaZa-17,Ur-17,Ur-17,ZhDuDu-19} for the singularity \eqref{g-sing}.

Looking in the above-mentioned contributions and, in general, in the literature, we notice that the natural question of uniqueness of solutions has received very few attention. More precisely, in the framework of indefinite problems, the uniqueness of positive solution is proved in \cite{BaPoTe-88,BrHe-90} when dealing with a concave $g(u)$, and in \cite{BrTo-10} when dealing with a singularity of the form \eqref{g-sing} with $\kappa=3$.
An intermediate situation is studied in \cite{HaZa-16ampa}; other types of special nonlinearities (convex-concave) are analysed in \cite{Na-16}.
We highlight that all these results concerning the uniqueness of the solutions are obtained for the linear differential operator.

As is shown by \cite{CoMa-89,ErTa-97,ErTa-98}, for convex or superlinear nonlinearities, the problem of uniqueness of positive solutions is of greater complexity even when the weight function $a(t)$ is of constant positive sign and apparently it has not yet been completely solved for sign-changing weights.
Indeed, for weight functions with multiple changes of the sign, uniqueness is not possible, in view of the results about the multiplicity of positive solutions obtained in
\cite{Fe-18book,FeZa-17jde}.

The aim of this paper is twofold: on one side, we plan to produce a uniqueness result including both linear and nonlinear differential operators; on the other hand, we investigate a situation allowing the superlinear nonlinearities as a special case.
Due to the fact that a weight function $a(t)$ with more than one change of sign allows multiplicity of positive solutions, it is natural to consider a coefficient $a(t)$ with a single change of sign. Similar coefficients have been considered in \cite{BoZa-15,BrTo-10,GoZa-19,LGOm-20,LGOmRi-17na,SoZa-17}. In particular, following \cite{BrTo-10,GoZa-19}, we will focus our attention on a stepwise weight function of the form
\begin{equation}\label{hp-a}
a(t) =
\begin{cases}
\, a_{+}, &\text{if $t\in\mathopen{[}0,\tau\mathclose{[}$,} \\
\, -a_{-}, &\text{if $t\in\mathopen{[}\tau,T\mathclose{[}$,}
\end{cases}
\end{equation}
where $a_{+},a_{-}>0$ and $0<\tau<T$. This framework allows us to study the uniqueness question exploiting techniques typical of autonomous system.

The general statement will be given in Section~\ref{section-2} (cf.~Theorem~\ref{th-phi}); by now, we just present it for some special cases. When dealing with the linear differential operator, i.e.~$\phi(s)=s$, it yields the following.

\begin{theorem}\label{th-u''}
Let $a\in L^{\infty}(0,T)$ be as in \eqref{hp-a}. Let $\gamma\in\mathopen{]}-\infty,-3\mathclose{]}\cup\mathopen{]}1,+\infty\mathclose{[}$. Then, the Neumann and the periodic boundary value problems associated with equation
\begin{equation}\label{eq-u''}
u'' + a(t) u^{\gamma}  = 0
\end{equation}
have at most one positive solution. Moreover, there exists a unique positive solution if and only if $\gamma\cdot\int_{0}^{T}a(t)\,\mathrm{d}t<0$.
\end{theorem}

Notice that the case of a singularity with $\gamma=-3$, already solved in \cite{BrTo-10}, is included in the above result. We mention that with our strategy of proof we can also deal with the linear case (i.e.~$\gamma=1$) so as to recover the existence of a simple principal eigenvalue (see Remark~\ref{rem-3.1} and Remark~\ref{rem-3.2}).

As for the more general case of a $p$-Laplacian operator $\phi(s)=\varphi_{p}(s)=|s|^{p-2}s$ with $p>1$, our main contribution is the following.

\begin{theorem}\label{th-p}
Let $a\in L^{\infty}(0,T)$ be as in \eqref{hp-a}. Let $p>1$ and
\begin{equation*}
\gamma\in\biggl{]}-\infty,\frac{1-2p}{p-1}\biggr{]} \cup \biggl{]}p-1,+\infty\biggr{[}.
\end{equation*}
Then, the Neumann and the periodic boundary value problems associated with equation
\begin{equation}\label{eq-p}
\bigl{(} |u'|^{p-2} u' \bigr{)}' + a(t) u^{\gamma}  = 0
\end{equation}
have at most one positive solution. Moreover, there exists a unique positive solution if and only if $\gamma\cdot\int_{0}^{T}a(t)\,\mathrm{d}t<0$.
\end{theorem}

The plan of the paper is the following. In Section~\ref{section-2}, we present our main abstract uniqueness result for equation~\eqref{eq-phi} and, to prove Theorem~\ref{th-u''} and Theorem~\ref{th-p}, in Section~\ref{section-3} we apply it to the case $\phi(s)=s$ and $\phi(s)=\varphi_{p}(s)=|s|^{p-2}s$ with $p>1$. In Section~\ref{section-4}, some remarks and open questions are presented, including a brief discussion for the Minkowski-curvature operator.

\section{An abstract uniqueness result}\label{section-2}

In this section, we aim to present a method to deal with a general class of nonlinear differential problems. Accordingly, we deal with the second-order equation
\begin{equation}\label{eq-main}
\bigl{(} \phi(u') \bigr{)}' + a(t) g(u) = 0,
\end{equation}
where, for $\Omega\subseteq\mathbb{R}$ an open interval with $0\in\Omega$, we assume that
\begin{itemize}
\item $\phi\colon\Omega\to\phi(\Omega)=\mathbb{R}$ is a homeomorphism with $\phi(0)=0$, $\phi(s)s>0$ for all $s\in\Omega\setminus\{0\}$;
\item $a\colon\mathopen{[}0,T\mathclose{[}\to\mathbb{R}$ is a step-wise function of the form \eqref{hp-a};
\item $g\colon\mathopen{]}0,+\infty\mathclose{[}\to\mathopen{]}0,+\infty\mathclose{[}$ is a continuous function.
\end{itemize}
In Theorem~\ref{th-phi} we state a uniqueness result for the Neumann boundary value problem associated with \eqref{eq-main}; a variant is given in Theorem~\ref{th-phi-bis}. In Section~\ref{section-2.2}, we discuss the periodic boundary value problem.

Preliminarily, we recall that a positive solution to equation~\eqref{eq-main} is a function $u\colon\mathopen{[}0,T\mathclose{]}\to\mathopen{]}0,+\infty\mathclose{[}$ of class $\mathcal{C}^{1}$ such that $\phi(u')$ is an absolutely continuous function and equation~\eqref{eq-main} is satisfied for
almost every $t\in\mathopen{[}0,T\mathclose{]}$.

\begin{remark}\label{rem-2.1}
When dealing with Neumann and periodic boundary conditions and with $g\in\mathcal{C}^{1}(\mathopen{]}0,+\infty\mathclose{[})$, an integration by parts of equation~\eqref{eq-main} divided by $g(u)$ gives
\begin{align*}
\bar{a} &:= a_{+}\tau-a_{-}(T-\tau)
= \int_{0}^{T} a(t) \,\mathrm{d}t
= -\int_{0}^{T} \dfrac{\bigl{(}\phi(u'(t))\bigr{)}'}{g(u(t))} \,\mathrm{d}t
\\ &\;= -\biggl{[} \dfrac{\phi(u'(t))}{g(u(t))} \biggr{]}_{0}^{T} - \int_{0}^{T} \dfrac{\phi(u'(t))u'(t)g'(u(t))}{(g(u(t)))^{2}} \,\mathrm{d}t.
\end{align*}
As a consequence, we immediately obtain a necessary condition for the existence of positive solutions of~\eqref{eq-main}, that is $\bar{a}<0$ if $g'(u)>0$ for all $u\in\mathopen{]}0,+\infty\mathclose{[}$, and $\bar{a}>0$ if $g'(u)<0$ for all $u\in\mathopen{]}0,+\infty\mathclose{[}$.
\hfill$\lhd$
\end{remark}

\subsection{The Neumann problem}\label{section-2.1}

Let us consider the planar system associated with equation \eqref{eq-main}, that is
\begin{equation}\label{syst}
\begin{cases}
\, x' = \phi^{-1}(y), \\
\, y' = -a(t) g(x).
\end{cases}
\end{equation}
A solution of \eqref{syst} is a couple $(x,y)$ of absolutely continuous functions satisfying \eqref{syst} for almost every $t$. Throughout the section, we confine ourselves in the half-right part $\mathopen{]}0,+\infty\mathclose{[}\times\mathbb{R}$ of the phase-plane.

According to the assumptions on $\phi(s)$, $a(t)$ and $g(x)$, for every time $t_{0}\in\mathopen{[}0,T\mathclose{[}$ and every initial condition $(x_{0},y_{0})\in\mathopen{]}0,+\infty\mathclose{[}\times\mathbb{R}$, system \eqref{syst} admits a unique local non-continuable solution with $x(t_{0})=x_{0}$ and $y(t_{0})=y_{0}$, denoted by
\begin{equation*}
\bigl{(} x(t;t_{0},x_{0},y_{0}),y(t;t_{0},x_{0},y_{0}) \bigr{)}.
\end{equation*}
The uniqueness of the solutions of the Cauchy problems is guaranteed by the special choice of the step-wise coefficient $a(t)$, indeed we enter the setting of the result in \cite{Re-00} concerning planar Hamiltonian systems. Moreover, we remark that $x(t;t_{0},x_{0},y_{0})>0$ for all $t$ in the maximal interval of existence, where the solution is defined.

In this section, we focus our attention on the Neumann boundary value problem associated with system~\eqref{syst}. Our goal is to prove that
\begin{quote}
\textit{there exists a unique initial condition $(x_{*},y_{*})\in\mathopen{]}0,+\infty\mathclose{[}\times\mathbb{R}$ such that $(x(t;\tau,x_{*},y_{*}),y(t;\tau,x_{*},y_{*}))$ is a solution of system~\eqref{syst} defined in $\mathopen{[}0,T\mathclose{]}$, $x(t;\tau,x_{*},y_{*})>0$ for all $t\in\mathopen{[}0,T\mathclose{]}$, and
\begin{equation}\label{N-bc}
y(0;\tau,x_{*},y_{*})=0, \qquad
y(T;\tau,x_{*},y_{*})=0.
\end{equation}}
\end{quote}

Setting $h=\phi^{-1}$, we can exploit hypothesis \eqref{hp-a} and so define the two autonomous systems
\begin{equation*}
(\mathcal{S}_{+})
\qquad
\begin{cases}
\, x' = h(y), \\
\, y' = - a_{+} g(x),
\end{cases}
\qquad
\qquad
(\mathcal{S}_{-})
\qquad
\begin{cases}
\, x' = h(y), \\
\, y' = a_{-} g(x).
\end{cases}
\end{equation*}

Let us assume that there is an hypothetical solution
\begin{equation*}
(x(t),y(t))=(x(t;\tau,x_{*},y_{*}),y(t;\tau,x_{*},y_{*}))
\end{equation*}
of~\eqref{syst} defined in $\mathopen{[}0,T\mathclose{]}$, with $x(t)>0$ for all $t\in\mathopen{[}0,T\mathclose{]}$, and satisfying the boundary conditions \eqref{N-bc}.
For future convenience, we set
\begin{equation}\label{eq-alphabeta}
\alpha = x(0), \qquad
\beta = x(T).
\end{equation}
Clearly, $\alpha>0$ and $\beta>0$. From $(\mathcal{S}_{+})$ and $y(0)=0$, we deduce that
\begin{equation*}
y(t) = y(0) + \int_{0}^{t} y'(\xi)\,\mathrm{d}\xi = - \int_{0}^{t} a_{+}g(x(\xi))\,\mathrm{d}\xi <0, \quad \text{for every $t\in\mathopen{]}0,\tau\mathclose{]}$,}
\end{equation*}
and thus $x'(t)=h(y(t))<0$ for every $t\in\mathopen{]}0,\tau\mathclose{[}$. Hence, $x(t)$ is strictly monotone decreasing in $\mathopen{[}0,\tau\mathclose{]}$ and $x(t)<\alpha$ for every $t\in\mathopen{]}0,\tau\mathclose{]}$. Analogously, from $(\mathcal{S}_{-})$ and $y(T)=0$, we have that $y(t)<0$ for every $t\in\mathopen{[}\tau,T\mathclose{[}$, so $x(t)$ is strictly monotone decreasing in $\mathopen{[}\tau,T\mathclose{]}$ and so $x(t)>\beta$ for every $t\in\mathopen{[}\tau,T\mathclose{[}$. We conclude that
\begin{equation}\label{y-neg}
y(t) < 0, \quad \text{for all $t\in\mathopen{]}0,T\mathclose{[}$,}
\end{equation}
so the solution $(x(t),y(t))$ is in the fourth quadrant $\mathopen{]}0,+\infty\mathclose{[}\times\mathopen{]}-\infty,0\mathclose{]}$. Moreover, $x(t)$ is strictly monotone decreasing in $\mathopen{[}0,T\mathclose{]}$, and thus
\begin{equation}\label{order}
0 < \beta < x_{*} < \alpha.
\end{equation}

Let $H$ and $G$ be primitives of $h$ and $g$, respectively. For convenience, we suppose that $H(0)=0$. When $g$ can be continuously extended to zero, we shall also assume $G(0)=0$.
We denote $G_{0} :=\lim_{x\to0}G(x)$ and $G_{+\infty} :=\lim_{x\to+\infty}G(x)$. In general, for an arbitrary function $g\colon\mathopen{]}0,+\infty\mathclose{[}\to\mathopen{]}0,+\infty\mathclose{[}$ the following four cases are possible:
\begin{itemize}
\item[$(i)$] $G_{0}\in\mathbb{R}$ and $G_{+\infty}=+\infty$;
in this situation, without loss of generality, we can suppose that $G\colon\mathopen{]}0,+\infty\mathclose{[}\to\mathopen{]}0,+\infty\mathclose{[}$ is a strictly monotone increasing surjective function;
\item[$(ii)$] $G_{0}\in\mathbb{R}$ and $G_{+\infty}\in\mathbb{R}$;
in this situation, without loss of generality, we can suppose that $G\colon\mathopen{]}0,+\infty\mathclose{[}\to\mathopen{]}0,L\mathclose{[}$ is a strictly monotone increasing surjective function, with $L:=\int_{0}^{+\infty}g(u)\,\mathrm{d}u$;
\item[$(iii)$]  $G_{0}=-\infty$ and $G_{+\infty}=+\infty$;
in this situation, without loss of generality, we can suppose that $G\colon\mathopen{]}0,+\infty\mathclose{[}\to\mathbb{R}$ is a strictly monotone increasing surjective function;
\item[$(iv)$] $G_{0}=-\infty$ and $G_{+\infty}\in\mathbb{R}$;
in this situation, without loss of generality, we can suppose that $G\colon\mathopen{]}0,+\infty\mathclose{[}\to\mathopen{]}-\infty,0\mathclose{[}$ is a strictly monotone increasing surjective function.
\end{itemize}

We observe that the quantities
$H(y) + a_{+} G(x)$ and $H(y) - a_{-} G(x)$ are constant for all $(x,y)$ solving
$(\mathcal{S}_{+})$ and $(\mathcal{S}_{-})$, respectively.
In particular, due to \eqref{N-bc}, \eqref{eq-alphabeta} and $H(0)=0$, we have that the solution $(x(t),y(t))$ satisfies
\begin{equation}\label{eq-levels}
\begin{aligned}
&H(y) + a_{+} G(x) = a_{+} G(\alpha), & &\text{on $\mathopen{[}0,\tau\mathclose{[}$,}&
\\
&H(y) - a_{-} G(x) = - a_{-} G(\beta), & &\text{on $\mathopen{[}\tau,T\mathclose{[}$.}&
\end{aligned}
\end{equation}

We notice that the functions $H_{l}:=H|_{\mathopen{]}-\infty,0\mathclose{]}}$,  $H_{r}:=H|_{\mathopen{[}0,+\infty\mathclose{[}}$ and $G$ are invertible since strictly monotone. For sake of simplicity in the notation, we set
\begin{equation*}
\mathcal{L}_{h} = h \circ H_{l}^{-1},
\qquad
\mathcal{L}_{g} = g \circ G^{-1}.
\end{equation*}
We remark that, when $h$ is odd, we find that $H_{l}^{-1}=-H_{r}^{-1}$ and, therefore, $\mathcal{L}_{h}=-h \circ H_{r}^{-1}$.

Recalling \eqref{y-neg}, from \eqref{eq-levels}, $(\mathcal{S}_{+})$ and $(\mathcal{S}_{-})$, we infer that
\begin{equation*}
\begin{aligned}
&x' = h(y) = \mathcal{L}_{h} ( a_{+} G(\alpha) - a_{+} G(x) ), & &\text{on $\mathopen{[}0,\tau\mathclose{[}$,}
\\
&x' = h(y) = \mathcal{L}_{h} ( a_{-} G(x) - a_{-} G(\beta) ), & &\text{on $\mathopen{[}\tau,T\mathclose{[}$.}
\end{aligned}
\end{equation*}
Recalling \eqref{order}, by an integration in $\mathopen{[}0,\tau\mathclose{[}$ and in $\mathopen{[}\tau,T\mathclose{[}$, we obtain
\begin{equation}\label{eq-star}
\begin{aligned}
&\int_{x_{*}}^{\alpha} \dfrac{\mathrm{d}x}{-\mathcal{L}_{h} ( a_{+} G(\alpha) - a_{+} G(x) )} = \tau,
\\
&\int_{\beta}^{x_{*}} \dfrac{\mathrm{d}x}{-\mathcal{L}_{h} ( a_{-} G(x) - a_{-} G(\beta) )} = T-\tau,
\end{aligned}
\end{equation}
respectively.
Conversely, if $\alpha$, $\beta$, $x_{*}$, with $0<\beta<x_{*}<\alpha$, are such that the above relations hold, we infer that the solution $(x(t;\tau,x_{*},y_{*}),y(t;\tau,x_{*},y_{*}))$ of the Cauchy problem is defined in $\mathopen{[}0,T\mathclose{]}$ and is such that $(x(0;\tau,x_{*},y_{*}),y(0;\tau,x_{*},y_{*}))=(\alpha,0)$ and $(x(T;\tau,x_{*},y_{*}),y(T;\tau,x_{*},y_{*}))=(\beta,0)$. Hence, by the choice of $h=\phi^{-1}$, we deduce that $x(t)$ is a positive decreasing solution of the Neumann boundary value problem associated with \eqref{eq-main}.

Let us perform the change of variable
\begin{equation*}
\vartheta = G(x), \quad x=G^{-1}(\vartheta), \quad
\mathrm{d}x = \dfrac{\mathrm{d}\vartheta}{g(G^{-1}(\vartheta))}=\dfrac{\mathrm{d}\vartheta}{\mathcal{L}_{g}(\vartheta)},
\end{equation*}
in \eqref{eq-star}, obtaining
\begin{equation}\label{eq-int+}
\int_{G(x_{*})}^{G(\alpha)} \dfrac{\mathrm{d}\vartheta}{-\mathcal{L}_{h} ( a_{+} G(\alpha) - a_{+} \vartheta ) \mathcal{L}_{g}(\vartheta)} = \tau
\end{equation}
and
\begin{equation}\label{eq-int-}
\int_{G(\beta)}^{G(x_{*})} \dfrac{\mathrm{d}\vartheta}{- \mathcal{L}_{h} ( a_{-} \vartheta - a_{-} G(\beta) ) \mathcal{L}_{g}(\vartheta)} = T-\tau.
\end{equation}
From \eqref{eq-levels} with $(x,y)=(x_{*},y_{*})$ one can deduce that
\begin{equation*}
G(x_{*}) = \dfrac{a_{+}G(\alpha)+a_{-}G(\beta)}{a_{+}+a_{-}},
\end{equation*}
thus $G(x_{*})$ is a convex combination of $G(\alpha)$ and $G(\beta)$.

Next, we set
\begin{equation*}
\omega := G(\alpha), \qquad
\sigma := G(\beta), \qquad
\mu := \dfrac{a_{+}}{a_{+}+a_{-}}.
\end{equation*}
Accordingly, \eqref{eq-int+} and \eqref{eq-int-} take the simplified form
\begin{equation}\label{eq-int+-bis}
\int_{\mu \omega+(1-\mu) \sigma}^{\omega} \dfrac{\mathrm{d}\vartheta}{-\mathcal{L}_{h} ( a_{+} \omega - a_{+} \vartheta ) \mathcal{L}_{g}(\vartheta)} = \tau
\end{equation}
and
\begin{equation}\label{eq-int--bis}
\int_{\sigma}^{\mu \omega+(1-\mu) \sigma} \dfrac{\mathrm{d}\vartheta}{- \mathcal{L}_{h} ( a_{-} \vartheta - a_{-} \sigma ) \mathcal{L}_{g}(\vartheta)} = T-\tau.
\end{equation}
Let $\mathcal{O}:=\{(\omega,\sigma)\in G(\mathopen{]}0,+\infty\mathclose{[})\times G(\mathopen{]}0,+\infty\mathclose{[}) \colon \omega > \sigma \}$ and we introduce the functions $\mathcal{M}_{\mathrm{I}},\mathcal{M}_{\mathrm{II}}\colon \mathcal{O}\to\mathopen{]}0,+\infty\mathclose{[}$ defined as
\begin{align*}
&\mathcal{M}_{\mathrm{I}}(\omega,\sigma):= \int_{\mu \omega+(1-\mu) \sigma}^{\omega} \dfrac{\mathrm{d}\vartheta}{-\mathcal{L}_{h} ( a_{+} \omega - a_{+} \vartheta ) \mathcal{L}_{g}(\vartheta)},
\\
&\mathcal{M}_{\mathrm{II}}(\omega,\sigma):= \int_{\sigma}^{\mu \omega+(1-\mu) \sigma} \dfrac{\mathrm{d}\vartheta}{- \mathcal{L}_{h} ( a_{-} \vartheta - a_{-} \sigma ) \mathcal{L}_{g}(\vartheta)}.
\end{align*}

From the above discussion, the following result holds true.

\begin{theorem}\label{th-phi}
Let $\Omega\subseteq\mathbb{R}$ be an open interval with $0\in\Omega$. Let  $\phi\colon\Omega \to\phi(\Omega)=\mathbb{R}$ be a homeomorphism with $\phi(0)=0$, $\phi(s)s>0$ for all $s\in\Omega\setminus\{0\}$. Let $a\in L^{\infty}(0,T)$ be as in \eqref{hp-a}. Let $g\colon\mathopen{]}0,+\infty\mathclose{[}\to\mathopen{]}0,+\infty\mathclose{[}$ be a continuous function.
Then, there exists a pair $(\omega,\sigma)\in\mathcal{O}$ which solves
\begin{equation}\label{syst-I-II-M}
\begin{cases}
\, \mathcal{M}_{\mathrm{I}}(\omega,\sigma)=\tau,
\\
\, \mathcal{M}_{\mathrm{II}}(\omega,\sigma)=T-\tau,
\end{cases}
\end{equation}
if and only if the Neumann boundary value problem associated with equation \eqref{eq-main} has a positive solution. Moreover, the unique solvability of \eqref{syst-I-II-M} is equivalent to the uniqueness of the positive solution of the Neumann problem.
\end{theorem}

\medskip

For our applications, we will consider a simplified but equivalent formulation of system~\eqref{syst-I-II-M} which can be obtained when $\omega=G(\alpha)$ is of constant sign. This excludes only the case $(iii)$ in the list above.

We introduce the new variable
\begin{equation*}
\rho := \dfrac{a_{-}G(\beta)}{a_{+}G(\alpha)} = \dfrac{a_{-}\sigma}{a_{+}\omega}.
\end{equation*}
Thus, we have
\begin{equation*}
\sigma = \dfrac{a_{+}}{a_{-}} \omega \rho, \qquad  G(x_{*}) = \mu \omega+(1-\mu) \sigma = \frac{a_{+}}{a_{+}+a_{-}} \, \omega \, (\rho+1).
\end{equation*}
By performing the change of variable $\vartheta= G(\alpha) \xi = \omega \xi$, formulas \eqref{eq-int+-bis} and \eqref{eq-int--bis} read as
\begin{align*}
& \omega \int_{\frac{a_{+}}{a_{+}+a_{-}} (\rho+1)}^{1} \dfrac{\mathrm{d}\xi}{-\mathcal{L}_{h} ( a_{+} \omega (1-\xi) ) \mathcal{L}_{g}(\omega \xi)} = \tau,
\\
& \omega \int_{\frac{a_{+}}{a_{-}} \rho}^{\frac{a_{+}}{a_{+}+a_{-}} (\rho+1)} \dfrac{\mathrm{d}\xi}{- \mathcal{L}_{h} ( \omega (a_{-} \xi - a_{+} \rho) ) \mathcal{L}_{g}(\omega \xi)} = T-\tau,
\end{align*}
respectively.

Let
\begin{equation}\label{domainD}
\mathcal{D}:=
\begin{cases}
\, G(\mathopen{]}0,+\infty\mathclose{[})\times\biggl{]}0,\dfrac{a_{-}}{a_{+}}\biggr{[},
&\text{in case $(i)$ and $(ii)$,}
\vspace{4pt}
\\
\, G(\mathopen{]}0,+\infty\mathclose{[})\times\biggl{]}\dfrac{a_{-}}{a_{+}},+\infty\biggr{[},
&\text{in case $(iv)$,}
\end{cases}
\end{equation}
and define the functions $\mathcal{F}_{\mathrm{I}},\mathcal{F}_{\mathrm{II}}\colon\mathcal{D}\to \mathopen{]}0,+\infty\mathclose{[}$ as follows
\begin{align*}
&\mathcal{F}_{\mathrm{I}}(\omega,\rho):= \omega \int_{\frac{a_{+}}{a_{+}+a_{-}} (\rho+1)}^{1} \dfrac{\mathrm{d}\xi}{-\mathcal{L}_{h} ( a_{+} \omega (1-\xi) ) \mathcal{L}_{g}(\omega \xi)},
\\
&\mathcal{F}_{\mathrm{II}}(\omega,\rho):= \omega \int_{\frac{a_{+}}{a_{-}} \rho}^{\frac{a_{+}}{a_{+}+a_{-}} (\rho+1)} \dfrac{\mathrm{d}\xi}{- \mathcal{L}_{h} ( \omega (a_{-} \xi - a_{+} \rho) ) \mathcal{L}_{g}(\omega \xi)}.
\end{align*}
From the above discussion, we have the following uniqueness result.

\begin{corollary}\label{cor-phi}
Let $\Omega\subseteq\mathbb{R}$ be an open interval with $0\in\Omega$. Let  $\phi\colon\Omega \to\phi(\Omega)=\mathbb{R}$ be a homeomorphism with $\phi(0)=0$, $\phi(s)s>0$ for all $s\in\Omega\setminus\{0\}$. Let $a\in L^{\infty}(0,T)$ be as in \eqref{hp-a}. Let $g\colon\mathopen{]}0,+\infty\mathclose{[}\to\mathopen{]}0,+\infty\mathclose{[}$ be a continuous function such that $(iii)$ does not hold.
Then, there exists a pair $(\omega,\rho)\in\mathcal{D}$ which solves
\begin{equation}\label{syst-I-II}
\begin{cases}
\, \mathcal{F}_{\mathrm{I}}(\omega,\rho)=\tau,
\\
\, \mathcal{F}_{\mathrm{II}}(\omega,\rho)=T-\tau,
\end{cases}
\end{equation}
if and only if the Neumann boundary value problem associated with equation \eqref{eq-main} has a positive solution. Moreover, the unique solvability of \eqref{syst-I-II} is equivalent to the uniqueness of the positive solution of the Neumann problem.
\end{corollary}

\medskip

We conclude this section, by presenting an equivalent version of Theorem~\ref{th-phi}.
Instead of \eqref{eq-levels}, our new starting point is the fact that the hypothetical solution $(x(t),y(t))$ also satisfies
\begin{equation*}
\begin{aligned}
&H(y) + a_{+} G(x) = H(y_{*}) + a_{+} G(x_{*}), & &\text{on $\mathopen{[}0,\tau\mathclose{[}$,}&
\\
&H(y) - a_{-} G(x) = H(y_{*}) - a_{-} G(x_{*}), & &\text{on $\mathopen{[}\tau,T\mathclose{[}$.}&
\end{aligned}
\end{equation*}
We immediately infer that
\begin{equation*}
\begin{aligned}
&y' = - a_{+}g(x) = - a_{+} \mathcal{L}_{g} \biggl{(} G(x_{*}) +\dfrac{H(y_{*})}{a_{+}} - \dfrac{H(y)}{a_{+}}  \biggr{)}, & &\text{on $\mathopen{[}0,\tau\mathclose{[}$,}
\\
&y' = a_{-}g(x) = a_{-} \mathcal{L}_{g} \biggl{(} G(x_{*}) +\dfrac{H(y)}{a_{-}} - \dfrac{H(y_{*})}{a_{-}}  \biggr{)}, & &\text{on $\mathopen{[}\tau,T\mathclose{[}$,}
\end{aligned}
\end{equation*}
and integrations in $\mathopen{[}0,\tau\mathclose{[}$ and, respectively, in $\mathopen{[}\tau,T\mathclose{[}$ give
\begin{equation}\label{syst-I-II-bis}
\begin{cases}
\, \displaystyle \int_{y_{*}}^{0} \dfrac{\mathrm{d}y}{a_{+} \mathcal{L}_{g} \biggl{(} G(x_{*}) +\dfrac{H(y_{*})}{a_{+}} - \dfrac{H(y)}{a_{+}}  \biggr{)}}=\tau,
\\
\, \displaystyle \int_{y_{*}}^{0} \dfrac{\mathrm{d}y}{a_{-} \mathcal{L}_{g} \biggl{(} G(x_{*}) +\dfrac{H(y)}{a_{-}} - \dfrac{H(y_{*})}{a_{-}}  \biggr{)}}=T-\tau.
\end{cases}
\end{equation}
The following result holds true.

\begin{theorem}\label{th-phi-bis}
Let $\Omega\subseteq\mathbb{R}$ be an open interval with $0\in\Omega$. Let  $\phi\colon\Omega \to\phi(\Omega)=\mathbb{R}$ be a homeomorphism with $\phi(0)=0$, $\phi(s)s>0$ for all $s\in\Omega\setminus\{0\}$. Let $a\in L^{\infty}(0,T)$ be as in \eqref{hp-a}. Let $g\colon\mathopen{]}0,+\infty\mathclose{[}\to\mathopen{]}0,+\infty\mathclose{[}$ be a continuous function.
Then, there exists a pair $(x_{*},y_{*})\in\mathopen{]}0,+\infty\mathclose{[}\times\mathopen{]}-\infty,0\mathclose{[}$ which solves \eqref{syst-I-II-bis}, if and only if the Neumann boundary value problem associated with equation \eqref{eq-main} has a positive solution. Moreover, the unique solvability of \eqref{syst-I-II-bis} is equivalent to the uniqueness of the positive solution of the Neumann problem.
\end{theorem}

\subsection{The periodic problem}\label{section-2.2}

In this section, we deal with the periodic boundary value problem associated with \eqref{eq-main} and we show that Theorem~\ref{th-phi} holds true also in the periodic case.
Following a procedure which is standard in this situation, we extend by $T$-periodicity the weigh $a(t)$ as an $L^{\infty}$-stepwise function defined in the whole real line.
In this framework, finding a solution of \eqref{eq-main} satisfying $u(0)=u(T)$ and $u'(0)=u'(T)$ is equivalent to finding a $T$-periodic solution of \eqref{eq-main} defined on $\mathbb{R}$.

As in Section~\ref{section-2.1}, we analyse the associated planar system \eqref{syst} and we look for periodic solutions
$(x(t),y(t))$ of \eqref{syst} such that $x(t)>0$ for all $t\in\mathbb{R}$.
Our purpose is to reduce the study of the periodic problem to the Neumann one, analysed in previous section. To this aim, we further assume that
\begin{center}
$h=\phi^{-1} \colon\mathbb{R}\to\mathbb{R}$ is odd.
\end{center}

The next two claims relate the existence/uniqueness of positive solutions of the $T$-periodic problem to the corresponding one for the Neumann problem.

\smallskip
\noindent\textit{Claim~1. Let $(x(t),y(t))$ be a solution of \eqref{syst} defined in the interval $\mathopen{[}\frac{\tau}{2},\frac{T+\tau}{2}\mathclose{]}$ with $x(t)>0$ for all $t\in\mathopen{[}\frac{\tau}{2},\frac{T+\tau}{2}\mathclose{]}$ and satisfying the Neumann condition at the boundary, that is
\begin{equation}\label{eq-y-N}
y\biggl{(}\dfrac{\tau}{2} \biggr{)} = y \biggl{(}\dfrac{T+\tau}{2}\biggl{)} = 0.
\end{equation}
Let also $(\hat{x}(t),\hat{y}(t))$ be the $T$-periodic extension of $(x(t),y(t))$ symmetric with respect to $t=\tau/2$, namely
\begin{equation}\label{extension}
(\hat{x}(t),\hat{y}(t))=
\begin{cases}
\, (x(t),y(t)), & \text{if $t\in\biggl{[}\dfrac{\tau}{2},\dfrac{T+\tau}{2}\biggr{]}$,}\vspace{2pt}
\\
\, (x(\tau-t),-y(\tau-t)), & \text{if $t\in\biggl{[}\dfrac{T-\tau}{2},\dfrac{\tau}{2}\biggr{]}$.}
\end{cases}
\end{equation}
Then, $(\hat{x}(t),\hat{y}(t))$ is a $T$-periodic solution of \eqref{syst} with $\hat{x}(t)>0$ for all $t\in\mathbb{R}$.}

Indeed, by construction, $(\hat{x}(t),\hat{y}(t))$ is symmetric with respect to $\frac{\tau}{2}$ and, by direct inspection, one can easily check that it is a solution of system \eqref{syst} on the interval $\mathopen{[}\frac{T-\tau}{2},\frac{T+\tau}{2}\mathclose{]}$; this follows from the fact that the extension of $a(t)$ by $T$-periodicity is symmetric with respect to $\frac{\tau}{2}$. Moreover, $(\hat{x}(t),\hat{y}(t))$ satisfies the $T$-periodic condition at the boundary of $\mathopen{[}\frac{T-\tau}{2},\frac{T+\tau}{2}\mathclose{]}$, that is
\begin{equation*}
\biggl{(}\hat{x}\biggl{(}\dfrac{\tau-T}{2}\biggr{)},\hat{y}\biggl{(}\dfrac{\tau-T}{2}\biggr{)}\biggr{)} = \biggl{(}\hat{x}\biggl{(}\dfrac{\tau+T}{2}\biggr{)},\hat{y}\biggl{(}\dfrac{\tau+T}{2}\biggr{)}\biggr{)}.
\end{equation*}
Since the weight $a(t)$ has been extended by $T$-periodicity on the whole real line and $(\hat{x}(t),\hat{y}(t))$ is a $T$-periodic extension of \eqref{extension}, we immediately conclude that $(\hat{x}(t),\hat{y}(t))$ solves \eqref{syst} and is such that $\hat{x}(t)>0$ for all $t\in\mathbb{R}$, and, by construction, $(\hat{x}(0),\hat{y}(0))=(\hat{x}(T),\hat{y}(T))$. Then, Claim~1 is proved.

\smallskip
\noindent\textit{Claim~2.
Let $(x(t),y(t))$ be a $T$-periodic solution of \eqref{syst} with $x(t)>0$ for all $t\in\mathbb{R}$. Then, the restriction of $(x(t),y(t))$ to the interval $\mathopen{[}\frac{\tau}{2},\frac{T+\tau}{2}\mathclose{]}$ is a solution of \eqref{syst} with $x(t)>0$ for all $t\in\mathopen{[}\frac{\tau}{2},\frac{T+\tau}{2}\mathclose{]}$ and satisfying the Neumann boundary condition \eqref{eq-y-N}.}

Indeed, for $(x(t),y(t))$ as in the assumption and by the special form \eqref{hp-a} of the weight function $a(t)$, we have that $y'(t)<0$ in $\mathopen{]}0,\tau\mathclose{[}$ and $y'(t)>0$ in $\mathopen{]}\tau,T\mathclose{[}$, and so $y(t)$ is strictly monotone decreasing in $\mathopen{[}0,\tau\mathclose{]}$ and strictly monotone increasing in $\mathopen{[}\tau,T\mathclose{]}$. Since $h$ is strictly monotone increasing, we also deduce that $x'(t)$ is strictly monotone decreasing in $\mathopen{[}0,\tau\mathclose{]}$ and strictly monotone increasing in $\mathopen{[}\tau,T\mathclose{]}$.

We notice that if $x'(0)=0$ then we have $x'(t)<0$ in $\mathopen{]}0,\tau\mathclose{[}$ and, since $x'(T)=h(y(T))=h(y(0))=x'(0)=0$, we also have $x'(t)<0$ in $\mathopen{]}\tau,T\mathclose{[}$; thus $x(t)$ is strictly monotone decreasing in $\mathopen{[}0,T\mathclose{]}$, a contradiction with $x(0)=x(T)$. Moreover, if $x'(\tau)=0$ then $x'(t)>0$ in $\mathopen{]}0,\tau\mathclose{[}$ and $x'(t)>0$ in $\mathopen{]}\tau,T\mathclose{[}$, a contradiction with $x(0)=x(T)$, as before. Therefore, from the above discussion and Rolle's theorem, we conclude that there exist exactly one critical point in $\mathopen{]}0,\tau\mathclose{[}$ and exactly one critical point in $\mathopen{]}\tau,T\mathclose{[}$. Let $\hat{t}\in\mathopen{]}0,\tau\mathclose{[}$ and $\check{t}\in\mathopen{]}\tau,T\mathclose{[}$ be such that $x'(\hat{t})=x'(\check{t})=0$. We claim that $\hat{t}$ and $\check{t}$ are the maximum point and, respectively, the minimum point of $x(t)$ in $\mathopen{[}0,T\mathclose{]}$. This is obvious since
$x(t)=x(\hat{t})+\int_{\hat{t}}^{t} x'(\xi) \,\mathrm{d}\xi < x(\hat{t})$ in every neighborhood of $\hat{t}$ (contained in $\mathopen{]}0,\tau\mathclose{[}$), while $x(t)=x(\check{t})+\int_{\check{t}}^{t} x'(\xi) \,\mathrm{d}\xi > x(\check{t})$ in every neighborhood of $\check{t}$ (contained in $\mathopen{]}\tau,T\mathclose{[}$).

We remark that the function $(x(t),y(t))$ is such that $y(\hat{t})=y(\check{t})=0$, $y(t)<0$ in $\mathopen{]}\hat{t},\check{t}\mathclose{[}$ and $x(t)$ is strictly monotone decreasing in $\mathopen{[}\hat{t},\check{t}\mathclose{]}$,  $y(t)>0$ in $\mathopen{]}0,\hat{t}\mathclose{[} \cup \mathopen{]}\check{t},T\mathclose{[}$ and $x(t)$ is strictly monotone increasing in $\mathopen{[}0,\hat{t}\mathclose{]} \cup \mathopen{[}\check{t},T\mathclose{]}$.

We are going to prove that $\hat{t}=\frac{\tau}{2}$ and $\check{t}=\frac{T+\tau}{2}$. Our approach is based on an analysis of the phase-portrait in the $(x,y)$-plane and is similar to the one performed in \cite{BrTo-10}. Due to the more general framework and in order to justify the additional hypothesis on $h$, we give here all the details.

Arguing as in Section~\ref{section-2.1} (cf.~\eqref{eq-levels}), we deduce that the solution $(x(t),y(t))$ satisfies
\begin{align}
&H(y) + a_{+} G(x) = a_{+} G(\alpha), \qquad \text{on $\mathopen{[}0,\tau\mathclose{]}$,}
\label{ll+}
\\
&H(y) - a_{-} G(x) = - a_{-} G(\beta), \qquad \text{on $\mathopen{[}\tau,T\mathclose{]}$,}
\label{ll-}
\end{align}
with $\alpha=x(\hat{t})$ and $\beta=x(\check{t})$.
The trajectory on the time interval $\mathopen{[}0,\hat{t}\mathclose{]}$ satisfies the relation
\begin{equation}\label{curve-1}
y = H_{r}^{-1}( a_{+} G(\alpha) - a_{+} G(x) ), \qquad \text{$(x,y)\in\mathopen{[}x(0),\alpha\mathclose{]}\times\mathopen{[}0,+\infty\mathclose{[}$,}
\end{equation}
which describes a strictly monotone decreasing curve. Analogously, using the fact that $H_{l}^{-1}=-H_{r}^{-1}$, the trajectory on the time interval $\mathopen{[}\hat{t},\tau\mathclose{]}$ satisfies the relation
\begin{equation}\label{curve-2}
y = - H_{r}^{-1}( a_{+} G(\alpha) - a_{+} G(x) ), \qquad \text{$(x,y)\in\mathopen{[}x(\tau),\alpha\mathclose{]}\times\mathopen{]}-\infty,0\mathclose{]}$,}
\end{equation}
which describes a strictly monotone increasing curve.
On the other hand, the trajectory on the time interval $\mathopen{[}\tau,\check{t}\mathclose{]}$ satisfies the relation
\begin{equation}\label{curve-3}
y = -H_{r}^{-1}( a_{-} G(x) - a_{-} G(\beta) ), \qquad \text{$(x,y)\in\mathopen{[}\beta,x(\tau)\mathclose{]}\times\mathopen{]}-\infty,0\mathclose{]}$,}
\end{equation}
which describes a strictly monotone decreasing curve, while the trajectory on the time interval $\mathopen{[}\check{t},T\mathclose{]}$ satisfies the relation
\begin{equation}\label{curve-4}
y = H_{r}^{-1}( a_{-} G(x) - a_{-} G(\beta) ), \qquad \text{$(x,y)\in\mathopen{[}\beta,x(T)\mathclose{]}\times\mathopen{[}0,+\infty\mathclose{[}$,}
\end{equation}
which describes a strictly monotone increasing curve.
By the strict monotonicity of the above curves, we infer that there is at most one intersection point between \eqref{curve-1} and \eqref{curve-4} in $\mathopen{]}\alpha,\beta\mathclose{[}\times\mathopen{]}0,+\infty\mathclose{[}$, and likewise at most one intersection point between \eqref{curve-2} and \eqref{curve-3} in $\mathopen{]}\alpha,\beta\mathclose{[}\times\mathopen{]}-\infty,0\mathclose{[}$. Actually these intersections exist and are uniquely determined as $(x_{+},y_{+})=(x(0),y(0))=(x(T),y(T))$ (in the first quadrant) and $(x_{-},y_{-})=(x(\tau),y(\tau))$ (in the fourth quadrant).
Due to the symmetry of $H$, it follows that $x_{+}=x_{-}$ and $y_{+}=-y_{-}$. Moreover, from the Abelian-type integrals representing the time-mappings, we find that the time necessary to connect (along the level line \eqref{ll+}) $(x_{+},y_{+})$ to $(\alpha,0)$ coincides with the time necessary to connect (along the same level line) $(\alpha,0)$ to $(x_{-},y_{-})$. Hence, $\hat{t}=\frac{\tau}{2}$. The same argument, with respect to the level line \eqref{ll-} shows that that the times necessary to connect $(x_{-},y_{-})$ to $(\beta,0)$ and $(\beta,0)$ to $(x_{+},y_{+})$ are equal. Then, $\hat{t}=\frac{\tau+T}{2}$. This concludes the proof of Claim~2.

\smallskip

From the above discussion, we deduce that we can reduce the problem of proving the existence and uniqueness of a positive solution of the $T$-periodic problem to the study of a Neumann problem defined in $\mathopen{[}\frac{\tau}{2},\frac{\tau+T}{2}\mathclose{]}$ for a step-wise function. Clearly, this latter problem is equivalent to the original problem studied in the previous section. In particular, observe also that systems \eqref{syst-I-II-M} and \eqref{syst-I-II} would be formally changed to corresponding new systems in which the target vector $(\tau,T-\tau)$ should be replaced by $(\frac{\tau}{2},\frac{T-\tau}{2})$. With reference to the result obtained in the subsequent section, no relevant point has to be changed due also to the elementary fact that the ration of the two components of the target vector $(\tau,T-\tau)$ or, respectively, $(\frac{\tau}{2},\frac{T-\tau}{2})$ remains unchanged (see the second equations in \eqref{syst-u''} and \eqref{syst-p}).

Finally, we conclude that Theorem~\ref{th-phi} and Corollary~\ref{cor-phi} are both valid also in the periodic case, with the additional assumption that $h$ is odd.

\section{Proofs of Theorem~\ref{th-u''} and Theorem~\ref{th-p}}\label{section-3}

In this section, we apply Corollary~\ref{cor-phi} when
\begin{equation}\label{def-g-2.3}
g(u)=u^{\gamma}, \quad u>0, \quad \gamma\in\mathbb{R},
\end{equation}
and for $\phi(s)=s$ in Section~\ref{section-3.1}, thus proving Theorem~\ref{th-u''}, while for $\phi(s)= |s|^{p-2}s$ with $p>1$ in Section~\ref{section-3.2}, thus proving Theorem~\ref{th-p}.
We focus our analysis on the Neumann boundary conditions; the result for the periodic problem follows from this as explained in Section~\ref{section-2.2} (indeed, the functions $h=\phi^{-1}$ will be odd).

\subsection{The case $\phi(s)=s$}\label{section-3.1}

We deal first with the simpler case
\begin{equation*}
\phi(s) = s.
\end{equation*}
As a first step, we prove that the Neumann boundary value problem associated with \eqref{eq-u''} has \textit{at most one} positive solution. As a second step, to conclude the proof, we show that there exists \textit{at least one} positive solution.

\medskip

In this special framework, we have
\begin{align*}
&h(y) = \phi^{-1}(y) = y, \qquad H(y) = \dfrac{y^{2}}{2}, \qquad y\in\mathbb{R},
\\
&H_{l}^{-1}(\xi) = -\sqrt{2\xi}, \qquad H_{r}^{-1}(\xi) = \sqrt{2\xi}, \qquad \xi\in\mathopen{[}0,+\infty\mathclose{[},
\end{align*}
which leads to
\begin{equation*}
\mathcal{L}_{h}(\xi) = -\sqrt{2\xi},
\qquad \xi\in\mathopen{[}0,+\infty\mathclose{[}.
\end{equation*}
Moreover, concerning the nonlinear term \eqref{def-g-2.3}, for $\gamma\in\mathbb{R}\setminus\{-1\}$, we deduce that
\begin{align}
&G(x) = \dfrac{x^{\gamma+1}}{\gamma+1}, \qquad x\in\mathopen{]}0,+\infty\mathclose{[},
\label{def-G-2.3}
\\
&G^{-1}(\xi) = |\gamma+1|^{\frac{1}{\gamma+1}} |\xi|^{\frac{1}{\gamma+1}},
\qquad \xi\in\mathrm{sign}(\gamma+1)\cdot\mathopen{]}0,+\infty\mathclose{[},
\label{def-G--2.3}
\end{align}
and thus
\begin{equation}\label{eq-GL}
\mathcal{L}_{g}(\xi) = |\gamma+1|^{\frac{\gamma}{\gamma+1}} |\xi|^{\frac{\gamma}{\gamma+1}},
\qquad \xi\in\mathrm{sign}(\gamma+1)\cdot\mathopen{]}0,+\infty\mathclose{[}.
\end{equation}
Incidentally, notice that in this situation the case $(iii)$ listed in Section~\ref{section-2.1} does not hold, so we are allowed to use Corollary~\ref{cor-phi}.

We study the unique solvability of system \eqref{syst-I-II} for $(\omega,\rho)\in\mathcal{D}$, where $\mathcal{D}$ is defined in \eqref{domainD}. In this special framework, we have
\begin{equation*}
\mathcal{F}_{\mathrm{I}}(\omega,\rho)
= \dfrac{|\gamma+1|^{-\frac{\gamma}{\gamma+1}}}{\sqrt{2}} \omega|\omega|^{-\frac{\gamma}{\gamma+1}-\frac{1}{2}}
\int_{\frac{a_{+}}{a_{+}+a_{-}}(\rho+1)}^{1} \dfrac{\mathrm{d}\xi}{ |a_{+} - a_{+}\xi |^{\frac{1}{2}} \, |\xi|^{\frac{\gamma}{\gamma+1}}}
\end{equation*}
and
\begin{align*}
\mathcal{F}_{\mathrm{II}}(\omega,\rho)
&= \dfrac{|\gamma+1|^{-\frac{\gamma}{\gamma+1}}}{\sqrt{2}}  \omega|\omega|^{-\frac{\gamma}{\gamma+1}-\frac{1}{2}}
\int_{\frac{a_{+}}{a_{-}}\rho}^{\frac{a_{+}}{a_{+}+a_{-}}(\rho+1)} \dfrac{\mathrm{d}\zeta}{ |a_{-}\zeta -a_{+}\rho |^{\frac{1}{2}} \, |\zeta|^{\frac{\gamma}{\gamma+1}} }
\\
&= \dfrac{|\gamma+1|^{-\frac{\gamma}{\gamma+1}}}{\sqrt{2}} \omega|\omega|^{-\frac{\gamma}{\gamma+1}-\frac{1}{2}} \rho^{1-\frac{\gamma}{\gamma+1}-\frac{1}{2}} \int_{\frac{a_{+}}{a_{-}}}^{\frac{a_{+}}{a_{+}+a_{-}}\frac{\rho+1}{\rho}}\dfrac{ \mathrm{d}\xi}{ |a_{-}\xi-a_{+} |^{\frac{1}{2}} \, |\xi|^{\frac{\gamma}{\gamma+1}}}.
\end{align*}

For simplicity in notation, we introduce the functions $I_{1},I_{2}\colon\mathopen{]}0,+\infty\mathclose{[}\to\mathbb{R}$ defined as
\begin{align*}
&I_{1}(\rho) = \int_{\frac{a_{+}}{a_{+}+a_{-}}(\rho+1)}^{1} \dfrac{\mathrm{d}\xi}{ |a_{+} - a_{+}\xi |^{\frac{1}{2}} \, |\xi|^{\frac{\gamma}{\gamma+1}}},
\\
&I_{2}(\rho)= \int_{\frac{a_{+}}{a_{-}}}^{\frac{a_{+}}{a_{+}+a_{-}}\frac{\rho+1}{\rho}} \dfrac{ \mathrm{d}\xi}{ |a_{-}\xi-a_{+} |^{\frac{1}{2}}\,|\xi|^{\frac{\gamma}{\gamma+1}}}.
\end{align*}
The functions $I_{1}$ and $I_{2}$ are formally defined for $\rho\in\mathopen{]}0,+\infty\mathclose{[}$, however in the applications we will study them on the intervals $\mathopen{]}0,a_{-}/a_{+}\mathclose{[}$ and $\mathopen{]}a_{-}/a_{+},+\infty\mathclose{[}$ depending on the choice of $\gamma$.
With this position, proving the unique solvability of system~\eqref{syst-I-II} for $(\omega,\rho)\in\mathcal{D}$ is equivalent to prove that there exists a unique pair $(\omega,\rho)\in\mathcal{D}$ which solves
\begin{equation}\label{syst-u''}
\begin{cases}
\,  \dfrac{|\gamma+1|^{-\frac{\gamma}{\gamma+1}}}{\sqrt{2}} \omega|\omega|^{-\frac{\gamma}{\gamma+1}-\frac{1}{2}} I_{1}(\rho) = \tau,
\\
\, \rho^{-1+\frac{\gamma}{\gamma+1}+\frac{1}{2}} \dfrac{I_{1}(\rho)}{I_{2}(\rho)} = \dfrac{\tau}{T-\tau},
\end{cases}
\end{equation}
where the second equation is the quotient between the two equation in~\eqref{syst-I-II}.

We now introduce the function
\begin{equation*}
F_{2}(\rho) := \rho^{\frac{\gamma-1}{2(\gamma+1)}} \dfrac{I_{1}(\rho)}{I_{2}(\rho)}
\end{equation*}
defined in $\mathopen{]}0,a_{-}/a_{+}\mathclose{[}$, for $\gamma>-1$, and in $\mathopen{]}a_{-}/a_{+},+\infty\mathclose{[}$, for $\gamma<-1$. We claim that $F_{2}(\rho)$ is strictly monotone
increasing for $\gamma\in\mathopen{]}-\infty,-3\mathclose{]}\cup\mathopen{[}1,+\infty\mathclose{[}$.

In the following computations, the prime symbol $'$ denotes the derivative with respect to $\rho$. We compute
\begin{align*}
I_{1}'(\rho) &= \dfrac{-\biggl{(} \dfrac{a_{+}}{a_{+}+a_{-}} (\rho+1) \biggr{)}' }{ \biggl{|} a_{+} - a_{+} \dfrac{a_{+}}{a_{+}+a_{-}} (\rho+1) \biggr{|}^{\frac{1}{2}} \, \biggl{(} \dfrac{a_{+}}{a_{+}+a_{-}}(\rho+1) \biggr{)}^{\!\frac{\gamma}{\gamma+1}}}
\\
&= - \biggl{(} \dfrac{a_{+}}{a_{+}+a_{-}} \biggr{)}^{\! 1-\frac{\gamma}{\gamma+1}-\frac{1}{2}} \dfrac{1}{ |a_{-}-a_{+}\rho|^{\frac{1}{2}} \, (\rho+1)^{\frac{\gamma}{\gamma+1}}  }
\end{align*}
and
\begin{align*}
I_{2}'(\rho)  &= \dfrac{\biggl{(} \dfrac{a_{+}}{a_{+}+a_{-}}\dfrac{\rho+1}{\rho} \biggr{)}' }{ \biggl{|} a_{-}\dfrac{a_{+}}{a_{+}+a_{-}}\dfrac{\rho+1}{\rho} - a_{+} \biggr{|}^{\frac{1}{2}} \biggl{(} \dfrac{a_{+}}{a_{+}+a_{-}}\dfrac{\rho+1}{\rho} \biggr{)}^{\!\frac{\gamma}{\gamma+1}}}
\\
&= - \rho^{-2+\frac{\gamma}{\gamma+1}+\frac{1}{2}} \biggl{(} \dfrac{a_{+}}{a_{+}+a_{-}} \biggr{)}^{\! 1-\frac{\gamma}{\gamma+1}-\frac{1}{2}} \dfrac{1}{ |a_{-} (\rho+1) - \rho(a_{+}+a_{-})|^{\frac{1}{2}} \, (\rho+1)^{\frac{\gamma}{\gamma+1}}}
\\
&= \rho^{-2+\frac{\gamma}{\gamma+1}+\frac{1}{2}} I_{1}'(\rho)
= \rho^{-\frac{\gamma+3}{2(\gamma+1)}} I_{1}'(\rho).
\end{align*}
Therefore,
\begin{equation*}
\dfrac{\mathrm{d}}{\mathrm{d}\rho} \biggl{(} \dfrac{I_{1}(\rho)}{I_{2}(\rho)} \biggr{)}
= \dfrac{I_{1}'(\rho)I_{2}(\rho)-I_{1}(\rho)I_{2}'(\rho)}{(I_{2}(\rho))^{2}}
= \dfrac{I_{1}'(\rho)}{(I_{2}(\rho))^{2}} \Bigl{(} I_{2}(\rho) - \rho^{-\frac{\gamma+3}{2(\gamma+1)}} I_{1}(\rho) \Bigr{)}.
\end{equation*}
We observe that
\begin{equation}\label{lim-0}
\lim_{\rho\to \frac{a_{-}}{a_{+}}} \Bigl{(} I_{2}(\rho) - \rho^{-\frac{\gamma+3}{2(\gamma+1)}} I_{1}(\rho) \Bigr{)} = 0
\end{equation}
and
\begin{align*}
\dfrac{\mathrm{d}}{\mathrm{d}\rho} \Bigl{(} I_{2}(\rho) - \rho^{-\frac{\gamma+3}{2(\gamma+1)}} I_{1}(\rho) \Bigr{)}
&= I_{2}'(\rho) + \dfrac{\gamma+3}{2(\gamma+1)} \rho^{-\frac{\gamma+3}{2(\gamma+1)}-1} I_{1}(\rho) - \rho^{-\frac{\gamma+3}{2(\gamma+1)}} I_{1}'(\rho)
\\
&= \dfrac{\gamma+3}{2(\gamma+1)} \rho^{-\frac{\gamma+3}{2(\gamma+1)}-1} I_{1}(\rho).
\end{align*}
We remark that $I_{1}(\rho)$ is defined in $\mathopen{]}0,a_{-}/a_{+}\mathclose{[}$ and is positive when $\gamma>-1$, while
$I_{1}(\rho)$ is defined in $\mathopen{]}a_{-}/a_{+},+\infty\mathclose{[}$ and is negative when $\gamma<-1$. Therefore, the derivative of the function $I_{2}(\rho) - \rho^{-\frac{\gamma+3}{2(\gamma+1)}}I_{1}(\rho)$ (in the interval of definition) has constant sign equal to the sign of $(\gamma+3)/(\gamma+1)$ on $\mathopen{]}0,a_{-}/a_{+}\mathclose{[}$ and to the sign of $-(\gamma+3)/(\gamma+1)$ on $\mathopen{]}a_{-}/a_{+},+\infty\mathclose{[}$. Using \eqref{lim-0}, we have that the function $I_{2}(\rho) - \rho^{-\frac{\gamma+3}{2(\gamma+1)}}I_{1}(\rho)$ (in the interval of definition) has constant sign equal to the sign of $-(\gamma+3)/(\gamma+1)$.
Finally, recalling that $I_{1}'(\rho)$ is negative, we infer that the (positive) function $I_{1}(\rho)/I_{2}(\rho)$ is strictly monotone increasing if $(\gamma+3)/(\gamma+1)>0$ and, respectively, is strictly monotone decreasing if $(\gamma+3)/(\gamma+1)<0$.

Summing up, since $F_{2}(\rho)$ is a product of two positive function, we have that $F_{2}(\rho)$ is strictly monotone if
\begin{equation*}
\dfrac{\gamma-1}{\gamma+1} \cdot \dfrac{\gamma+3}{\gamma+1} > 0,
\end{equation*}
thus, if
\begin{equation*}
\gamma < -3 \quad \text{or} \quad \gamma > 1.
\end{equation*}
We stress that in these cases $F_{2}(\rho)$ is strictly monotone increasing.
Moreover, if $\gamma=-3$ we obtain that the function $I_{1}(\rho)/I_{2}(\rho)$ is constant and so $F_{2}(\rho)=\rho$ is strictly monotone increasing, while if $\gamma=1$ then $F_{2}(\rho)=I_{1}(\rho)/I_{2}(\rho)$ is strictly monotone increasing.

From the above discussion we conclude that if $\gamma\in\mathopen{]}-\infty,-3\mathclose{]}\cup\mathopen{[}1,+\infty\mathclose{[}$, then the function $F_{2}(\rho)$ is strictly monotone increasing. The claim is proved.

Hence, we have shown that there exists at most one value of $\rho$ such that the second equation in \eqref{syst-u''} holds. When $\gamma\neq1$ the first equation in \eqref{syst-u''} is uniquely solvable, and thus for $\gamma\in\mathopen{]}-\infty,-3\mathclose{]}\cup\mathopen{]}1,+\infty\mathclose{[}$ there exists at most one pair $(\omega,\rho)$ which solves \eqref{syst-u''}.

To conclude the proof of Theorem~\ref{th-u''}, in view of Remark~\ref{rem-2.1}, we only need to prove that if $\gamma\cdot \bar{a}<0$, then there exists at least one pair $(\omega,\rho)$ which solves \eqref{syst-u''}. Preliminarily, we recall that if the second equation (which depends only on $\rho$) has a solution $\hat{\rho}$, then the first equation is solvable, since its first member (with $\rho=\hat{\rho}$) is surjective on $\mathopen{]}0,+\infty\mathclose{[}$.

We observe that the function $F_{2}(\rho)$ is continuous and strictly monotone increasing in $\mathopen{]}0,a_{-}/a_{+}\mathclose{[}$, for $\gamma> 1$, and, respectively, in $\mathopen{]}a_{-}/a_{+},+\infty\mathclose{[}$, for $\gamma\leq -3$. We compute the limits as $\rho$ approaches the endpoints. By L'H\^{o}pital's rule, we have
\begin{equation*}
\lim_{\rho\to \frac{a_{-}}{a_{+}}} \dfrac{I_{1}(\rho)}{I_{2}(\rho)}
= \lim_{\rho\to\frac{a_{-}}{a_{+}}} \dfrac{I_{1}'(\rho)}{I_{2}'(\rho)}
= \lim_{\rho\to \frac{a_{-}}{a_{+}}} \rho^{-\frac{\gamma-1}{2(\gamma+1)}+1}
= \biggl{(}\dfrac{a_{-}}{a_{+}}\biggr{)}^{\!-\frac{\gamma-1}{2(\gamma+1)}+1}
\end{equation*}
and thus
\begin{equation}\label{limit-a-+}
\lim_{\rho\to \frac{a_{-}}{a_{+}}} F_{2}(\rho) = \dfrac{a_{-}}{a_{+}}.
\end{equation}
If $\gamma> 1$, then
\begin{equation}\label{lim-00}
\lim_{\rho\to 0^{+}} F_{2}(\rho) = 0.
\end{equation}
Indeed, as $\rho\to0^{+}$, $I_{1}(\rho)$ and $I_{2}(\rho)$ tend to positive constants, and $\rho^{\frac{\gamma-1}{2(\gamma+1)}}\to0^{+}$. Therefore, we conclude that the second equation in \eqref{syst-u''} has a solution if and only if $\frac{a_{-}}{a_{+}} > \frac{\tau}{T-\tau}$ which is equivalent to assume that $\bar{a}<0$ (cf.~Remark~\ref{rem-2.1}).
While, if $\gamma\leq -3$, then
\begin{equation*}
\lim_{\rho\to +\infty} F_{2}(\rho) = +\infty.
\end{equation*}
Indeed, as $\rho\to+\infty$, $I_{1}(\rho)$ and $I_{2}(\rho)$ tend to negative constants, and $\rho^{\frac{\gamma-1}{2(\gamma+1)}}\to+\infty$.
Therefore, we conclude that the second equation in \eqref{syst-u''} has a solution if and only if $\frac{a_{-}}{a_{+}} < \frac{\tau}{T-\tau}$ which is equivalent to assume that $\bar{a}>0$.

The proof of Theorem~\ref{th-u''} is completed.
\hfill\qed

\begin{remark}[The case $\gamma=1$]\label{rem-3.1}
The proof of Theorem~\ref{th-u''} in particular states that, for $\gamma=1$,  the function $F_{2}(\rho)$ is strictly monotone increasing. Arguing as in the second part of the proof, one can also show that \eqref{limit-a-+} and \eqref{lim-00} are still valid, using L'H\^{o}pital's rule and, respectively, the fact that, when $\gamma=1$, as $\rho\to0^{+}$, $I_{1}(\rho)$ tends to a positive constant, while $I_{2}(\rho)\to+\infty$. We can conclude that the second equation in \eqref{syst-u''} is uniquely solved. Let $\hat{\rho}$ be the solution, then the first equation in \eqref{syst-u''} reads as $I_{1}(\hat{\rho})/2 = \tau$, which either is not solvable, or holds for all $\omega\in \mathopen{]}0,+\infty\mathclose{[}$.
This is in agreement with the well known fact that the eigenspace of Neumann/periodic one-signed solutions of the linear equation $u''+a(t)u=0$ has dimension less than or equal to $1$ (cf.~\cite{CoLe-55}).
\hfill$\lhd$
\end{remark}

\begin{remark}[The linear/nonlinear eigenvalue problem]\label{rem-3.2}
Let us consider the nonlinear eigenvalue problem
\begin{equation*}
u'' + \lambda a(t) u^{\gamma} = 0
\end{equation*}
together with the Neumann/periodic boundary conditions. In this case, in the above analysis, system~\eqref{syst-u''} should be modified to
\begin{equation*}
\begin{cases}
\,  \dfrac{|\gamma+1|^{-\frac{\gamma}{\gamma+1}}}{\sqrt{2}} \omega|\omega|^{-\frac{\gamma}{\gamma+1}-\frac{1}{2}} I_{1}(\rho) = \sqrt{\lambda} \tau,
\\
\, F_{2}(\rho) = \dfrac{\tau}{T-\tau}.
\end{cases}
\end{equation*}
The study of the second equation of this system is the same as above and, for $\gamma\in\mathopen{]}-\infty,-3\mathclose{]}\cup\mathopen{[}1,+\infty\mathclose{[}$, it provides a unique solution $\hat{\rho}$. By inserting this value in the first equation, for $\gamma\neq1$, we obtain the curve
\begin{equation*}
\omega = G(u(0)) = \Biggl{(} \frac{\sqrt{2} \tau}{|\gamma+1|^{-\frac{\gamma}{\gamma+1}} I_{1}(\hat{\rho})} \Biggr{)}^{\!\frac{2(\gamma+1)}{1-\gamma}} \lambda^{\frac{\gamma+1}{1-\gamma}}
\end{equation*}
and thus the classical bifurcation diagrams appearing in similar situations (cf.~\cite[Case~1, p.~446]{Li-82}).
On the other hand, if $\gamma=1$, we have
\begin{equation*}
\lambda=\lambda_{1} := \biggl{(}\dfrac{I_{1}(\hat{\rho})}{2 \tau} \biggr{)}^{\!2},
\end{equation*}
which ensures the existence of a unique positive principal eigenvalue consistently with the classical linear theory (cf.~\cite{Bo-14,BrLi-80,Ze-05} and, concerning the $p$-Laplacian, see \cite{MeYaZh-10} and the references therein).
\hfill$\lhd$
\end{remark}

\subsection{The case $\phi(s)= |s|^{p-2}s$ with $p>1$}\label{section-3.2}

We deal with the more general case
\begin{equation*}
\phi(s)=\varphi_{p}(s) = |s|^{p-2}s, \quad p>1.
\end{equation*}
The proof of Theorem~\ref{th-p} follows similar steps as the proof of Theorem~\ref{th-u''} in Section~\ref{section-3.1}; for this reason, in order to avoid unnecessary repetitions, we focus only on the main differences and skip the verification of some computations.

Recalling that the inverse of $\varphi_{p}$ is $\varphi_{q}(s)=|s|^{q-2}s$, where $\frac{1}{p}+\frac{1}{q}=1$, we have
\begin{align*}
&h(y) = \varphi_{q}(y) = |y|^{\frac{2-p}{p-1}}y,
\qquad H(y) = \dfrac{p-1}{p}|y|^{\frac{p}{p-1}-1}y, \qquad y\in\mathbb{R},
\\
&H_{l}^{-1}(\xi) = - \biggl{(}\dfrac{p}{p-1}\biggr{)}^{\!\frac{p-1}{p}} |\xi|^{\frac{p-1}{p}},
\qquad \xi\in\mathopen{]}-\infty,0\mathclose{]},
\\
&H_{r}^{-1}(\xi) = \biggl{(}\dfrac{p}{p-1}\biggr{)}^{\!\frac{p-1}{p}} \xi^{\frac{p-1}{p}},
\qquad \xi\in\mathopen{[}0,+\infty\mathclose{[},
\end{align*}
which leads to
\begin{equation*}
\mathcal{L}_{h}(\xi) = - \biggl{(}\dfrac{p}{p-1}\biggr{)}^{\!\frac{1}{p}} |\xi|^{\frac{1}{p}},
\qquad \xi\in\mathopen{[}0,+\infty\mathclose{[}.
\end{equation*}
The functions $g(x)$, $G(x)$ and $G^{-1}(\xi)$ are define as in \eqref{def-g-2.3}, \eqref{def-G-2.3}, \eqref{def-G--2.3}, respectively. Hence, $\mathcal{L}_{g}$ is defined as in \eqref{eq-GL}.

With these positions, proceeding similarly as in Section~\ref{section-3.1}, we introduce the functions
\begin{align*}
&I_{1}(\rho) = \int_{\frac{a_{+}}{a_{+}+a_{-}}(\rho+1)}^{1} \dfrac{\mathrm{d}\xi}{ |a_{+} - a_{+}\xi |^{\frac{1}{p}} \, |\xi|^{\frac{\gamma}{\gamma+1}}},
\\
&I_{2}(\rho)= \int_{\frac{a_{+}}{a_{-}}}^{\frac{a_{+}}{a_{+}+a_{-}}\frac{\rho+1}{\rho}} \dfrac{ \mathrm{d}\xi}{ |a_{-}\xi-a_{+} |^{\frac{1}{p}}\,|\xi|^{\frac{\gamma}{\gamma+1}}},
\end{align*}
and, by Corollary~\ref{cor-phi}, we infer that proving that system~\eqref{syst-I-II} admits a unique solution $(\omega,\rho)\in\mathcal{D}$ is equivalent to prove that there exists a unique pair $(\omega,\rho)\in\mathcal{D}$ which solves
\begin{equation}\label{syst-p}
\begin{cases}
\,  \dfrac{|\gamma+1|^{-\frac{\gamma}{\gamma+1}}}{\biggl{(}\dfrac{p}{p-1}\biggr{)}^{\!\frac{1}{p}}} \omega|\omega|^{-\frac{\gamma}{\gamma+1}-\frac{1}{p}} I_{1}(\rho) = \tau,
\\
\, \rho^{-1+\frac{\gamma}{\gamma+1}+\frac{1}{p}} \dfrac{I_{1}(\rho)}{I_{2}(\rho)} = \dfrac{\tau}{T-\tau},
\end{cases}
\end{equation}
which is exactly system \eqref{syst-u''} when $p=2$.

As above, we study the auxiliary function
\begin{equation*}
F_{p}(\rho) = \rho^{-1+\frac{\gamma}{\gamma+1}+\frac{1}{p}} \dfrac{I_{1}(\rho)}{I_{2}(\rho)}.
\end{equation*}
We compute
\begin{align*}
&I_{1}'(\rho) = \biggl{(} \dfrac{a_{+}}{a_{+}+a_{-}} \biggr{)}^{1-\frac{\gamma}{\gamma+1}-\frac{1}{p}} \dfrac{1}{ |a_{-}-a_{+}\rho|^{\frac{1}{p}} \, (\rho+1)^{\frac{\gamma}{\gamma+1}} },
\\
&I_{2}'(\rho) = \rho^{-2+\frac{\gamma}{\gamma+1}+\frac{1}{p}} I_{1}'(\rho),
\end{align*}
and
\begin{equation*}
\dfrac{\mathrm{d}}{\mathrm{d}\rho} \biggl{(} \dfrac{I_{1}(\rho)}{I_{2}(\rho)} \biggr{)}
= \dfrac{I_{1}'(\rho)}{(I_{2}(\rho))^{2}} \Bigl{(} I_{2}(\rho) - \rho^{-2+\frac{\gamma}{\gamma+1}+\frac{1}{p}} I_{1}(\rho) \Bigr{)}.
\end{equation*}
Moreover, we have
\begin{align*}
&\lim_{\rho\to \frac{a_{-}}{a_{+}}} \Bigl{(}  I_{2}(\rho) - \rho^{-2+\frac{\gamma}{\gamma+1}+\frac{1}{p}} I_{1}(\rho) \Bigr{)} = 0,
\\
&\dfrac{\mathrm{d}}{\mathrm{d}\rho} \Bigl{(}  I_{2}(\rho) - \rho^{-2+\frac{\gamma}{\gamma+1}+\frac{1}{p}} I_{1}(\rho) \Bigr{)}
= \biggl{(}2-\frac{\gamma}{\gamma+1}-\frac{1}{p}\biggr{)} \rho^{-3+\frac{\gamma}{\gamma+1}+\frac{1}{p}} I_{1}(\rho).
\end{align*}
Reasoning as in the proof of Theorem~\ref{th-u''}, we conclude that $F_{p}(\rho)$ is strictly monotone if
\begin{equation*}
\biggl{(}-1+\frac{\gamma}{\gamma+1}+\frac{1}{p}\biggr{)}  \cdot \biggl{(}2-\frac{\gamma}{\gamma+1}-\frac{1}{p}\biggr{)}> 0,
\end{equation*}
which is equivalent to
\begin{equation*}
(p-1) \biggl{(}\gamma - \frac{1-2p}{p-1} \biggl{)} \bigl{(} \gamma - (p-1) \bigr{)} > 0.
\end{equation*}
Moreover, if $\gamma=(1-2p)/(p-1)$ we obtain that the function $I_{1}(\rho)/I_{2}(\rho)$ is constant and so $F_{p}$ is strictly monotone increasing, while if $\gamma=p-1$ then $F_{p}(\rho)=I_{1}(\rho)/I_{2}(\rho)$ is strictly monotone increasing.

For the part concerning the existence of solutions, we can complete the proof of Theorem~\ref{th-p} by arguing exactly as in Section~\ref{section-3.1}.
\hfill\qed

\begin{remark}\label{rem-3.3}
In this context the case $\gamma=p-1$ plays a similar role as the linear case $u''+a(t)u=0$. Similar considerations as those in Remark~\ref{rem-3.1} are in order.
\hfill$\lhd$
\end{remark}

\section{Final remarks}\label{section-4}

In this final section, we collect some complementary results and open questions motivated by our approach.

\subsection{The ``excluded'' values of $\gamma$}\label{section-4.1}

Even if Theorem~\ref{th-u''} and Theorem~\ref{th-p} cover a wide range of values of the real exponent $\gamma$, some interesting cases are not investigated.
In this section, we present some comments concerning these situations. In order to simplify the exposition, we focus our attention on the linear differential operator, namely $\phi(s)=s$; similar analysis could be performed for the $p$-Laplacian operator.

It is convenient to split the analysis for $\gamma\in\mathopen{]}-3,1\mathclose{[}\setminus\{-1\}$ into the two cases $\gamma\in\mathopen{]}-1,1\mathclose{[}$ and $\gamma\in\mathopen{]}-3,-1\mathclose{[}$. In the first situation, there may be non-negative solutions which vanish (when $0<\gamma<1$) or hit the singularity in finite time (when $-1<\gamma<0$) in the interval $\mathopen{[}\tau,T\mathclose{[}$ where the weight is negative. In the second case, the singularity satisfies the strong force condition (cf.~\cite{To-15,Ur-16}) and this prevents the possibility that solutions reach the singularity.

\smallskip
\noindent
\textit{Case 1. $\gamma\in\mathopen{]}-1,1\mathclose{[}$.}
Recalling system~\eqref{syst-u''}, we consider the function $F_{2}(\rho)$ with $\rho\in\mathopen{]}0,a_{-}/a_{+}\mathclose{[}$. Preliminarily, we notice that for $\gamma=0$ an easy computation shows that $F_{2}(\rho)\equiv a_{-}/a_{+}$ for all $\rho\in\mathopen{]}0,a_{-}/a_{+}\mathclose{[}$. Thus, from now on, we exclude this trivial situation.

By L'H\^{o}pital's rule, we observe that
\begin{align*}
\lim_{\rho\to0^{+}} \dfrac{\rho^{\frac{\gamma-1}{2(\gamma+1)}}}{I_{2}(\rho)} &= \lim_{\rho\to0^{+}} \dfrac{\gamma-1}{2(\gamma+1)} \dfrac{\rho^{-\frac{\gamma+3}{2(\gamma+1)}}}{I_{2}'(\rho)}
= \dfrac{\gamma-1}{2(\gamma+1)} \lim_{\rho\to0^{+}} \dfrac{\rho^{-\frac{\gamma+3}{2(\gamma+1)}}}{\rho^{-\frac{\gamma+3}{2(\gamma+1)}} I_{1}'(\rho)}
\\
&= \dfrac{\gamma-1}{2(\gamma+1)} \lim_{\rho\to0^{+}} \dfrac{1}{I_{1}'(\rho)}
= \dfrac{1-\gamma}{2(\gamma+1)} \biggl{(} \dfrac{a_{+}+a_{-}}{a_{+}} \biggr{)}^{\! 1-\frac{\gamma}{\gamma+1}-\frac{1}{2}} (a_{-})^{\frac{1}{2}}.
\end{align*}
Therefore, we deduce that
\begin{equation*}
\lim_{\rho\to0^{+}} F_{2}(\rho) =
\dfrac{1-\gamma}{2(\gamma+1)} \biggl{(} \dfrac{a_{+}+a_{-}}{a_{+}} \biggr{)}^{\! \frac{1-\gamma}{2(\gamma+1)}} \biggl{(}\dfrac{a_{-}}{a_{+}}\biggr{)}^{\!\frac{1}{2}}
\int_{\frac{a_{+}}{a_{+}+a_{-}}}^{1} \dfrac{\mathrm{d}\xi}{ |1 - \xi |^{\frac{1}{2}} \, |\xi|^{\frac{\gamma}{\gamma+1}}}=:K_{0}(\gamma).
\end{equation*}
On the other hand, the limit of $F_{2}(\rho)$ as $\rho\to (a_{-}/a_{+})^{-}$ is the same as obtained in \eqref{limit-a-+}.

The fact that $K_{0}(\gamma)>0$ for all $\gamma\in\mathopen{]}-1,1\mathclose{[}$ implies that the function $F_{2}(\rho)$ is not surjective on $\mathopen{]}0,a_{-}/a_{+}\mathclose{[}$. As a consequence, the condition $\gamma\cdot\bar{a}<0$ is not sufficient for the existence of positive solutions.
This fact is not surprising and actually it is consistent with the previous observation that in the interval where the weight is negative we have the possibility of solutions vanishing in finite time, due to the lack of uniqueness for the Cauchy problems, or hitting the singularity, due to the absence of the strong force condition in zero.

For further convenience, we claim that the function $K_{0}\colon\mathopen{]}-1,1\mathclose{[}\to\mathopen{]}0,+\infty\mathclose{[}$ is a strictly monotone decreasing function such that
\begin{equation}\label{K0}
\lim_{\gamma\to(-1)^{+}} K_{0}(\gamma) = +\infty, \quad
\lim_{\gamma\to1^{-}} K_{0}(\gamma) = 0, \quad
K_{0}(0)=\frac{a_{-}}{a_{+}}.
\end{equation}
Indeed, equalities in \eqref{K0} can be easily checked by a direct inspection. Thus, it is sufficient to check the strict monotonicity of the map $\psi\colon\mathopen{]}-1,1\mathclose{[}\to\mathopen{]}0,+\infty\mathclose{[}$ defined as
\begin{equation*}
\psi(\gamma) = \dfrac{1-\gamma}{2(\gamma+1)} L^{\frac{1-\gamma}{2(\gamma+1)}}
\int_{\frac{1}{L}}^{1} \dfrac{\mathrm{d}\xi}{ (1 - \xi)^{\frac{1}{2}} \, \xi^{\frac{\gamma}{\gamma+1}}},
\end{equation*}
where $L = (a_{+}+a_{-})/a_{-}>1$. Setting $\vartheta=(1-\gamma) / (2(\gamma+1))$, which maps $\mathopen{]}-1,1\mathclose{[}$ to $\mathopen{]}0,+\infty\mathclose{[}$ as a strictly monotone decreasing function, we are led to prove that the function $\eta\colon\mathopen{]}0,+\infty\mathclose{[}\to\mathopen{]}0,+\infty\mathclose{[}$ defined as
\begin{equation*}
\eta(\vartheta) = \vartheta \int_{\frac{1}{L}}^{1} \dfrac{(L\xi)^{\vartheta}}{ (1 - \xi )^{\frac{1}{2}} \, \xi^{\frac{1}{2}}}\,\mathrm{d}\xi
\end{equation*}
is strictly monotone increasing. Observing that the map $\vartheta \mapsto (L\xi)^{\vartheta}$ is strictly monotone increasing for every fixed $\xi\in\mathopen{]}1/L,1\mathclose{[}$, the claim is thus proved.

We are now ready to discuss the solvability of the boundary value problems.
In the case $\gamma\in\mathopen{]}0,1\mathclose{[}$, the nonlinearity is concave and smooth in $\mathopen{]}0,+\infty\mathclose{[}$ and we can apply \cite[Lemma~3.1]{BaPoTe-88} (see also \cite{BrHe-90}) to ensure the fact that there exists at most one positive solution of the problem. This implies that the function $F_{2}(\rho)$ is strictly monotone. More precisely, taking into account the monotonicity of $K_{0}(\gamma)$ and $K_{0}(0)= a_{-}/a_{+}$, we have
\begin{equation*}
\lim_{\rho\to0^{+}} F_{2}(\rho) = K_{0}(\gamma) < \dfrac{a_{-}}{a_{+}} = \lim_{\rho\to\frac{a_{-}}{a_{+}}} F_{2}(\rho).
\end{equation*}
Therefore, $F_{2}(\rho)$ is strictly monotone increasing.
Hence, taking into account the continuity of $F_{2}(\rho)$ on $\mathopen{]}0,a_{-}/a_{+}\mathclose{[}$, we can state the following.

\begin{proposition}\label{prop-4.1}
Let $a\in L^{\infty}(0,T)$ be as in \eqref{hp-a}. Let $\gamma\in \mathopen{]}0,1\mathclose{[}$. Then, the Neumann and the periodic boundary value problems associated with equations \eqref{eq-u''} and \eqref{eq-p} have a positive solution if and only if
\begin{equation}\label{cond-4.1}
K_{0}(\gamma) < \dfrac{\tau}{T-\tau} < \dfrac{a_{-}}{a_{+}}.
\end{equation}
Moreover, the solution is unique.
\end{proposition}

In the case $\gamma\in \mathopen{]}-1,0\mathclose{[}$, by the previous argument as before, we have instead
\begin{equation*}
\lim_{\rho\to\frac{a_{-}}{a_{+}}} F_{2}(\rho) = \dfrac{a_{-}}{a_{+}} <  K_{0}(\gamma) = \lim_{\rho\to0^{+}} F_{2}(\rho).
\end{equation*}
Even without information of the monotonicity, the continuity of $F_{2}(\rho)$ implies that its range covers the open interval $\mathopen{]}a_{-}/a_{+},K_{0}(\gamma)\mathclose{[}$. Hence, we have the following.

\begin{proposition}\label{prop-4.2}
Let $a\in L^{\infty}(0,T)$ be as in \eqref{hp-a}. Let $\gamma\in \mathopen{]}-1,0\mathclose{[}$. Then, the Neumann and the periodic boundary value problems associated with equations \eqref{eq-u''} and \eqref{eq-p} have a positive solution if
\begin{equation}\label{cond-4.2}
\dfrac{a_{-}}{a_{+}} < \dfrac{\tau}{T-\tau} < K_{0}(\gamma).
\end{equation}
\end{proposition}

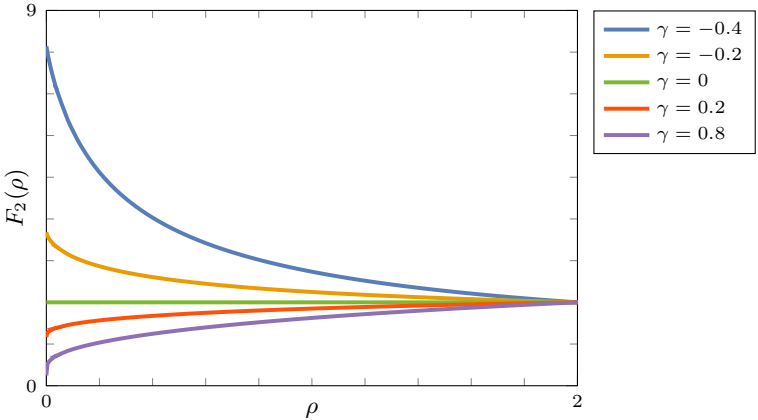
\begin{figure}[htb]
\centering
\begin{tikzpicture}[scale=1]
\begin{axis}[
  legend pos = outer north east,
  legend cell align={left},
  tick label style={font=\scriptsize},
          scale only axis,
  enlargelimits=false,
  xtick={0,2},
  xticklabels={$0$,$2$},
  ytick={0,9},
  yticklabels={$0$,$9$},
  xlabel={\small $\rho$},
  ylabel={\small $F_{2}(\rho)$},
  max space between ticks=50,
                minor x tick num=9,
                minor y tick num=8,
every axis x label/.style={
below,
at={(3.5cm,0cm)},
  yshift=-3pt
  },
every axis y label/.style={
below,
at={(0cm,2.5cm)},
  xshift=-3pt},
  y label style={rotate=90,anchor=south},
  width=7cm,
  height=5cm,
  xmin=0,
  xmax=2,
  ymin=0,
  ymax=9]
\addplot [draw=c1, line width=1.5pt, smooth] coordinates
{(0, 8.14276) (0.004, 7.99401) (0.008, 7.86468) (0.012, 7.74572) (0.016, 7.63453) (0.02, 7.52969) (0.024, 7.43027) (0.028, 7.33557) (0.032, 7.24508) (0.036, 7.15839) (0.04, 7.07515) (0.08, 6.38473) (0.12, 5.86655) (0.16, 5.45554) (0.2, 5.11821) (0.24, 4.83462) (0.28, 4.59182) (0.32, 4.38093) (0.36, 4.1956) (0.4, 4.03112) (0.44, 3.88393) (0.48, 3.75126) (0.52, 3.63092) (0.56, 3.52118) (0.6, 3.42061) (0.64, 3.32803) (0.68, 3.24247) (0.72, 3.16313) (0.76, 3.0893) (0.8, 3.0204) (0.84, 2.95593) (0.88, 2.89545) (0.92, 2.83858) (0.96, 2.78499) (1., 2.73439) (1.04, 2.68652) (1.08, 2.64116) (1.12, 2.5981) (1.16, 2.55717) (1.2, 2.5182) (1.24, 2.48104) (1.28, 2.44558) (1.32, 2.41168) (1.36, 2.37925) (1.4, 2.34818) (1.44, 2.31839) (1.48, 2.28979) (1.52, 2.26232) (1.56, 2.2359) (1.6, 2.21047) (1.64, 2.18598) (1.68, 2.16236) (1.72, 2.13958) (1.76, 2.11759) (1.8, 2.09634) (1.84, 2.0758) (1.88, 2.05592) (1.92, 2.03669) (1.96, 2.01806) (2,2)};
\addplot [draw=c2, line width=1.5pt, smooth] coordinates
{(0, 3.68105) (0.004, 3.59797) (0.008, 3.54896) (0.012, 3.50914) (0.016, 3.47468) (0.02, 3.4439) (0.024, 3.41589) (0.028, 3.39007) (0.032, 3.36603) (0.036, 3.34351) (0.04, 3.32228) (0.08, 3.15582) (0.12, 3.03658) (0.16, 2.94276) (0.2, 2.86529) (0.24, 2.79934) (0.28, 2.742) (0.32, 2.69136) (0.36, 2.64607) (0.4, 2.60518) (0.44, 2.56796) (0.48, 2.53384) (0.52, 2.50239) (0.56, 2.47325) (0.6, 2.44613) (0.64, 2.4208) (0.68, 2.39705) (0.72, 2.37472) (0.76, 2.35367) (0.8, 2.33377) (0.84, 2.31491) (0.88, 2.297) (0.92, 2.27996) (0.96, 2.26373) (1., 2.24823) (1.04, 2.23341) (1.08, 2.21922) (1.12, 2.20561) (1.16, 2.19255) (1.2, 2.18) (1.24, 2.16792) (1.28, 2.15628) (1.32, 2.14506) (1.36, 2.13424) (1.4, 2.12378) (1.44, 2.11367) (1.48, 2.10389) (1.52, 2.09442) (1.56, 2.08524) (1.6, 2.07634) (1.64, 2.06771) (1.68, 2.05932) (1.72, 2.05118) (1.76, 2.04326) (1.8, 2.03557) (1.84, 2.02807) (1.88, 2.02078) (1.92, 2.01367) (1.96, 2.00675) (2,2)};
\addplot [draw=c3, line width=1.5pt] coordinates
{(0,2) (2,2)};
\addplot [draw=c4, line width=1.5pt, smooth] coordinates
{(0, 1.16883) (0.004, 1.26484) (0.008, 1.29219) (0.012, 1.31175) (0.016, 1.32753) (0.02, 1.34099) (0.024, 1.35284) (0.028, 1.3635) (0.032, 1.37324) (0.036, 1.38224) (0.04, 1.39063) (0.08, 1.45458) (0.12, 1.49994) (0.16, 1.53612) (0.2, 1.56661) (0.24, 1.59314) (0.28, 1.61673) (0.32, 1.63802) (0.36, 1.65746) (0.4, 1.67537) (0.44, 1.692) (0.48, 1.70751) (0.52, 1.72207) (0.56, 1.73578) (0.6, 1.74875) (0.64, 1.76106) (0.68, 1.77276) (0.72, 1.78393) (0.76, 1.7946) (0.8, 1.80482) (0.84, 1.81463) (0.88, 1.82405) (0.92, 1.83313) (0.96, 1.84188) (1., 1.85032) (1.04, 1.85848) (1.08, 1.86638) (1.12, 1.87402) (1.16, 1.88144) (1.2, 1.88863) (1.24, 1.89561) (1.28, 1.9024) (1.32, 1.90901) (1.36, 1.91543) (1.4, 1.92169) (1.44, 1.92779) (1.48, 1.93374) (1.52, 1.93955) (1.56, 1.94522) (1.6, 1.95075) (1.64, 1.95616) (1.68, 1.96146) (1.72, 1.96663) (1.76, 1.9717) (1.8, 1.97666) (1.84, 1.98151) (1.88, 1.98627) (1.92, 1.99094) (1.96, 1.99551) (2,2)};
\addplot [draw=c5, line width=1.5pt, smooth] coordinates
{(0, 0.246849) (0.004, 0.5093) (0.008, 0.557925) (0.012, 0.59172) (0.016, 0.618628) (0.02, 0.641424) (0.024, 0.661438) (0.028, 0.679423) (0.032, 0.695849) (0.036, 0.711035) (0.04, 0.725204) (0.08, 0.834413) (0.12, 0.913885) (0.16, 0.978871) (0.2, 1.03488) (0.24, 1.08465) (0.28, 1.12973) (0.32, 1.17113) (0.36, 1.20955) (0.4, 1.24548) (0.44, 1.2793) (0.48, 1.31128) (0.52, 1.34166) (0.56, 1.37062) (0.6, 1.39831) (0.64, 1.42486) (0.68, 1.45037) (0.72, 1.47494) (0.76, 1.49864) (0.8, 1.52155) (0.84, 1.54372) (0.88, 1.5652) (0.92, 1.58604) (0.96, 1.60629) (1., 1.62598) (1.04, 1.64514) (1.08, 1.66381) (1.12, 1.68201) (1.16, 1.69977) (1.2, 1.71711) (1.24, 1.73405) (1.28, 1.75062) (1.32, 1.76683) (1.36, 1.7827) (1.4, 1.79823) (1.44, 1.81346) (1.48, 1.82839) (1.52, 1.84303) (1.56, 1.85739) (1.6, 1.8715) (1.64, 1.88535) (1.68, 1.89895) (1.72, 1.91232) (1.76, 1.92546) (1.8, 1.93839) (1.84, 1.95111) (1.88, 1.96362) (1.92, 1.97593) (1.96, 1.98806) (2,2)};
\legend{
\scriptsize{$\gamma = -0.4$}, \scriptsize{$\gamma = -0.2$}, \scriptsize{$\gamma = 0$}, \scriptsize{$\gamma = 0.2$}, \scriptsize{$\gamma = 0.8$}
}
\end{axis}
\end{tikzpicture}
\caption{Graphs of $F_{2}(\rho)$ for $\gamma\in\{-0.4,-0.2, 0, 0.2, 0.8\}$, with $a_{+}=1$ and $a_{-}=2$.}
\label{fig-01}
\end{figure}

To the best of our knowledge, there are no general uniqueness results available for this range of the exponent. Numerical simulations suggest that the function $F_{2}(\rho)$ is strictly monotone decreasing (see Figure~\ref{fig-01}). This would guarantee that the positive solution is unique also in this situation and consequently condition \eqref{cond-4.2} would be sharp, likewise condition \eqref{cond-4.1} is sharp in the case $\gamma\in \mathopen{]}0,1\mathclose{[}$.

To compare our result with a similar case previously studied, we recall that in \cite[Corollary~3]{GoZa-19jdde} the authors investigate a singular equation with a stepwise indefinite function, and obtain a range for the existence of positive solutions, that for instance in the case $\gamma=-1/2$ and $a_{-}=a_{+}$ would read like
\begin{equation*}
1 < \dfrac{\tau}{T-\tau} \leq \sqrt{\dfrac{35}{27}}.
\end{equation*}
On the other hand, an easy computation shows that  $K_{0}(\gamma)=5$ and thus \eqref{cond-4.2} reads as
\begin{equation*}
1 < \dfrac{\tau}{T-\tau} < 5.
\end{equation*}
This shows that, at least in this case, the range provided by our result is better.

\smallskip
\noindent
\textit{Case 2. $\gamma\in\mathopen{]}-3,-1\mathclose{[}$.}
We consider the function $F_{2}(\rho)$ with $\rho\in\mathopen{]}a_{-}/a_{+},+\infty\mathclose{[}$ and compute the limits at the endpoints. The limit of $F_{2}(\rho)$ as $\rho\to (a_{-}/a_{+})^{-}$ is the same as obtained in \eqref{limit-a-+}. As $\rho\to+\infty$, $I_{1}(\rho)$ and $I_{2}(\rho)$ tend to negative constants and, thus, $F_{2}(\rho)\to+\infty$ as $\rho\to+\infty$.
The continuity of $F_{2}(\rho)$ on $\mathopen{]}a_{-}/a_{+},+\infty\mathclose{[}$ implies that its range covers the open interval $\mathopen{]}a_{-}/a_{+},+\infty\mathclose{[}$. Therefore, we have the following.

\begin{proposition}\label{prop-4.3}
Let $a\in L^{\infty}(0,T)$ be as in \eqref{hp-a}. Let $\gamma\in \mathopen{]}-3,-1\mathclose{[}$. Then, the Neumann and the periodic boundary value problems associated with equations \eqref{eq-u''} and \eqref{eq-p} have a positive solution if and only if $\int_{0}^{T} a(t)\,\mathrm{d}t>0$.
\end{proposition}

\begin{figure}[htb]
\centering
\begin{tikzpicture}[scale=1]
\begin{axis}[
  legend pos = outer north east,
  legend cell align={left},
  tick label style={font=\scriptsize},
          scale only axis,
  enlargelimits=false,
  xtick={2,20},
  xticklabels={$2$,$20$},
  ytick={0,300},
  yticklabels={$0$,$300$},
  xlabel={\small $\rho$},
  ylabel={\small $F_{2}(\rho)$},
  max space between ticks=50,
                minor x tick num=17,
                minor y tick num=5,
every axis x label/.style={
below,
at={(3.5cm,0cm)},
  yshift=-3pt
  },
every axis y label/.style={
below,
at={(0cm,2.5cm)},
  xshift=-3pt},
  y label style={rotate=90,anchor=south},
  width=7cm,
  height=5cm,
  xmin=2,
  xmax=20,
  ymin=0,
  ymax=300]
\addplot [draw=c1, line width=1.5pt, smooth] coordinates
{(2, 2) (3., 4.5589) (4., 8.33606) (5., 13.4742) (6., 20.1049) (7., 28.3495) (8., 38.3214) (9., 50.127) (10., 63.8671) (11., 79.6372) (12., 97.5288) (13., 117.63) (14., 140.024) (15., 164.793) (16., 192.015) (17., 221.767) (18., 254.122) (19., 289.153) (20., 326.929)};
\addplot [draw=c2, line width=1.5pt, smooth] coordinates
{(2, 2) (3., 4.14974) (4., 7.05907) (5., 10.7525) (6., 15.2526) (7., 20.5801) (8., 26.7533) (9., 33.789) (10., 41.7025) (11., 50.5078) (12., 60.2181) (13., 70.8455) (14., 82.4016) (15., 94.8971) (16., 108.342) (17., 122.747) (18., 138.12) (19., 154.47) (20., 171.806)};
\addplot [draw=c3, line width=1.5pt, smooth] coordinates
{(2, 2) (3., 3.69256) (4., 5.74801) (5., 8.14263) (6., 10.8596) (7., 13.886) (8., 17.2111) (9., 20.8259) (10., 24.7229) (11., 28.895) (12., 33.3365) (13., 38.0418) (14., 43.0061) (15., 48.2249) (16., 53.694) (17., 59.4098) (18., 65.3687) (19., 71.5674) (20., 78.0029)};
\addplot [draw=c4, line width=1.5pt, smooth] coordinates
{(2, 2) (3., 3.44435) (4., 5.0882) (5., 6.90771) (6., 8.88642) (7., 11.012) (8., 13.2746) (9., 15.6662) (10., 18.1801) (11., 20.8105) (12., 23.5525) (13., 26.4016) (14., 29.3538) (15., 32.4058) (16., 35.5542) (17., 38.7962) (18., 42.1293) (19., 45.5508) (20., 49.0586)};
\addplot [draw=c5, line width=1.5pt, smooth] coordinates
{(2, 2) (20., 20.)};
\legend{
\scriptsize{$\gamma = -1.5$}, \scriptsize{$\gamma = -1.6$}, \scriptsize{$\gamma = -1.8$}, \scriptsize{$\gamma = -2$}, \scriptsize{$\gamma = -3$}
}
\end{axis}
\end{tikzpicture}
\caption{Graphs of $F_{2}(\rho)$ for $\gamma\in\{-3,-2,-1.8,-1.6,-1.5\}$, with $a_{+}=1$ and $a_{-}=2$.}
\label{fig-02}
\end{figure}
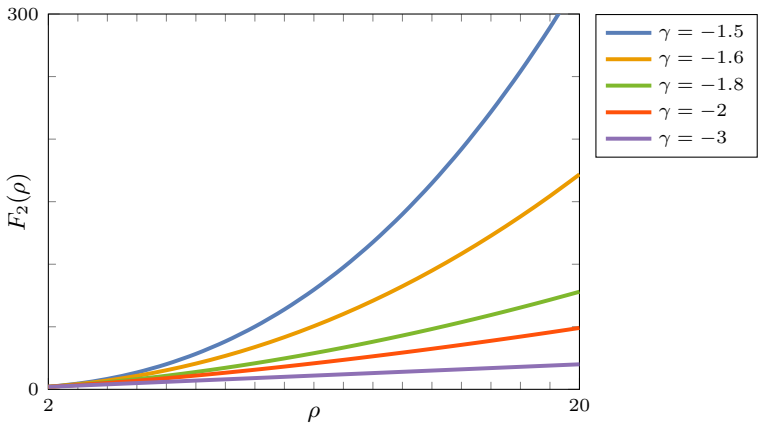

We observe that Proposition~\ref{prop-4.3} is already contained in \cite[Theorem~1]{Ur-16}. Our contribution is just in providing a different proof.
Concerning the problem of uniqueness, numerical simulations suggest that $F_{2}(\rho)$ is strictly monotone increasing in $\mathopen{]}a_{-}/a_{+},+\infty\mathclose{[}$ (see Figure~\ref{fig-02}). This would guarantee that the positive solution is unique.

At this point is natural to raise the following conjecture.

\begin{conjecture}\label{conj-1}
Let $a\in L^{\infty}(0,T)$ be as in \eqref{hp-a}. Let $\gamma\in\mathbb{R}\setminus\{-1,0,1\}$. Then, the Neumann and the periodic boundary value problems associated with equations \eqref{eq-u''} and \eqref{eq-p} have at most one positive solution.
\end{conjecture}

\subsection{The Minkowski-curvature operator}\label{section-4.2}

In this paper, we have applied Theorem~\ref{th-phi} to the case of homogeneous differential operators $\phi(s)$ and homogeneous nonlinearities $g(u)$. This choice has the advantage to simplify \eqref{syst-I-II-M} to \eqref{syst-u''} or \eqref{syst-p}, and the consequent analysis of the functions $F_{2}(\rho)$ or, respectively, $F_{p}(\rho)$.
However, our technique appears useful to study more general situations, especially when the function $\mathcal{L}_{h}$ can be computed explicitly. For example, in the case of the Minkowski-curvature equation
\begin{equation}\label{eq-Mink}
\Biggl{(} \dfrac{u'}{\sqrt{1-(u')^{2}}}\Biggr{)}' + a(t) g(u) = 0,
\end{equation}
we have
\begin{equation*}
\phi(s) = \dfrac{s}{\sqrt{1-s^{2}}}
\end{equation*}
and therefore, for $h=\phi^{-1}$, we find
\begin{equation*}
h(y) = \dfrac{y}{\sqrt{1+y^{2}}}, \quad H(y)=\sqrt{1+y^{2}}-1, \quad \mathcal{L}_{h}(\xi) = -\dfrac{\sqrt{\xi^{2}+2\xi}}{1+\xi}.
\end{equation*}
As a consequence, the analysis of system \eqref{syst-I-II-M} could be simplified for some special choices of $g(u)$. We do not pursue here this investigation which can be the topic of future researches.

Recent works for the Neumann and periodic problems associated with \eqref{eq-Mink} show that for $g(u)=u^{\gamma}$ with $\gamma>1$ multiple positive solutions do exist also in the case of a weight with a single change of sign as in \eqref{hp-a} (see \cite{BoFe-20na,BoFe-20jde}). On the other hand, numerical simulations suggest the possibility of uniqueness results when $g(u)$ is a strictly increasing function with ``super-exponential'' growth at infinity (cf.~\cite{BoFeZa-pp}).

\bibliographystyle{elsart-num-sort}
\bibliography{BoFeZa-biblio}

\end{document}